\documentclass{amsart}

\usepackage[margin=1.25in]{geometry}
\usepackage{amsmath,amssymb,amscd,mathtools}
\usepackage{booktabs,tabularx,array}
\usepackage{placeins}
\usepackage{microtype}
\usepackage[hidelinks]{hyperref}
\usepackage{tikz-cd}
\usepackage{quiver}

\allowdisplaybreaks[1]
\theoremstyle{plain}
\newtheorem{theorem}{Theorem}[section]
\newtheorem{proposition}[theorem]{Proposition}
\newtheorem{lemma}[theorem]{Lemma}
\newtheorem{corollary}[theorem]{Corollary}
\newtheorem{conjecture}[theorem]{Conjecture}
\theoremstyle{definition}
\newtheorem*{question}{Question}

\theoremstyle{remark}
\newtheorem{remark}[theorem]{Remark}

\numberwithin{equation}{section}

\newcommand{\cl}{\operatorname{cl}}
\newcommand{\expn}{\operatorname{exp}}
\newcommand{\M}{\mathrm M}
\newcommand{\B}{\widetilde{\mathrm B}_0}
\newcommand{\F}{\mathbb F}
\newcommand{\Z}{\mathbb Z}
\newcommand{\g}{\gamma}
\newcommand{\curlyw}{\mathbin{\curlywedge}}
\newcommand{\wed}{\mathbin{\wedge}}
\newcommand{\calP}{\mathcal P}
\newcommand{\Kset}{K}
\newcommand{\Aut}{\operatorname{Aut}}

\title[Exponent bounds for finite $p$-groups]%
{Schur multiplier exponents of finite
$p$-groups: bounds and critical-class phenomena}
\author[Z. Xing]{Zugan Xing}
\address{School of Mathematical Sciences, Shanghai Jiao Tong University\\
800 Dongchuan Road\\
Shanghai 200240, China}
\email{xingzugan@aliyun.com}

\author[S. Zheng]{Shaowu Zheng}
\address{Department of Mathematics, University of Auckland\\
Private Bag 92019\\
Auckland 1142, New Zealand}
\email{shaowu.zheng@auckland.ac.nz}

\hypersetup{
  pdftitle={Schur multiplier exponents of finite p-groups:
bounds and critical-class phenomena},
  pdfauthor={Zugan Xing and Shaowu Zheng},
  pdfkeywords={Schur multiplier, finite p-groups, exponent bounds, commutator subgroup, lower central series, Bogomolov multiplier}
}
\date{}

\begin{document}

\begin{abstract}
The Schur exponent conjecture is a long-standing conjecture in group theory, which asserts that \(
    \expn\M(G)\mid\expn G
\) for every finite group $G$, where $\M(G)$ denotes the Schur multiplier of $G$.

In this article, we prove that the Schur exponent conjecture holds for finite $p$-groups $G$ of nilpotency class at most \(2p-2\) when $p$ is an odd prime. Our approach is to establish upper bounds on the exponent of $\gamma_{s+1}(G)$ for $s\in\mathbb{N}$. Using this approach, we further establish exponent bounds for Schur multipliers, nonabelian exterior squares, and tensor squares of $G$ of arbitrary nilpotency class, improving some existing bounds. Moreover, we show that Moravec's Bogomolov multiplier exponent conjecture holds for finite \(p\)-groups of nilpotency class at most \(2p-1\).

Havas conjectured that, for every prime $p>3$, there exists a finite $p$-group of exponent $p$ and nilpotency class $2p-1$ for which the Schur exponent conjecture fails. We show that the Schur exponent conjecture still holds for finite $p$-groups of nilpotency class at most $2p-1$ generated by at most three elements. For four-generator groups of exponent $p$ and class at most $2p-1$, we prove that $\M(G)\cong C_{p^2}^{\,b}\times C_p^{\,a}$ for some $a\ge0$ and $b\in\{0,1\}$.
\end{abstract}

\maketitle

\section{Introduction}\label{sec:introduction}

The Schur multiplier $\M(G)$ of a finite group
$G$, originating from Schur's work on projective
representations \cite{Schur1904,Schur1907},
is a classical invariant in group theory.

In number theory, Schur multipliers play an important role in the study of Galois groups with restricted ramification and non-abelian Cohen--Lenstra heuristics 
\cite{BostonBushHajir2017,LiuXing2026Arxiv,LiuXing2026JNT,Watt1985}. 
Moreover, for a finite Galois extension $L/k$ with Galois group $G$, the exponent of \(M(G)\) provides an upper bound for the exponent of the Hasse norm knot group \(\mathfrak K(L/k)\) \cite[Corollary 3.7]{MacedoNewton2022}.

These applications motivate the study of the exponent of the Schur multiplier. A fundamental question is whether
\[
\expn\M(G)\mid\expn G.
\]
It is known to hold for powerful $p$-groups \cite{LubotzkyMann1987} and $p$-groups of maximal class \cite{Moravec2011}. Antony, Komma and Thomas proved that it holds for groups of nilpotency class at most $p$ \cite[Corollary~2.5]{AntonyKommaThomas2022}.  However, Vaughan-Lee constructed a $5$-group of class $9$ and exponent $5$ whose Schur multiplier has exponent $25$  \cite{VaughanLee2021counterexamples}. He further showed that if $G$ is a $5$-group of nilpotency class less than $9$, or a $7$-group of nilpotency class less than $13$, then \(\exp M(G)\mid \expn G\) \cite[Theorem 9]{VaughanLee2021II}. These results motivate the following question.

\begin{question}\label{ques:class-bound}
Let $p$ be an odd prime. Does \(
\expn\M(G)\mid\expn G
\) hold for every finite $p$-group $G$ with $\cl(G)\le 2p-2$? Moreover, is the upper bound $2p-2$ sharp? More precisely, does there exist, for every prime $p>3$, a finite $p$-group $G$ of nilpotency class $2p-1$ such that \(
    \expn\M(G)\nmid\expn G?
\)
\end{question}
We prove that the above divisibility holds. We further establish the following result.
\begin{theorem}\label{thm:exterior-tensor-bounds}
Let $p$ be an odd prime and let $G$ be a  finite nontrivial
$p$-group. Then:
\begin{enumerate}
\item If $\cl(G)\le 2p-2$, then $\expn\M(G)\mid\expn G$.
\item More generally, for $X\in\{G\wed G,\M(G),G\otimes G\}$, we have $
   \expn X\mid
   (\expn G)^{
      1+\left\lceil
      \log_p\frac{\cl(G)+1}{2p-1}
      \right\rceil}.
$
\end{enumerate}
\end{theorem}

The above result improves the exponent bounds for $G\otimes G$, $M(G)$, and $G\wed G$ previously obtained by Bastos, Melo, Gonçalves and Monetta \cite[Theorem~1.4]{BastosMeloGoncalvesMonetta2022}, Sambonet \cite[Theorem~1.1]{Sambonet2017}, and Antony, Komma and Thomas \cite[Theorem~4.2]{AntonyKommaThomas2022}, respectively. For a more detailed comparison, we refer to Table~\ref{tab:comparison-class-bounds}. 

The Bogomolov multiplier $\B(G)$ is another important invariant in group theory and number theory. Moravec conjectured \cite[Conjecture~1.3]{Moravec2019Bogomolov} that $\expn\B(G)\mid\expn G$ for every finite group $G$. 
In Section~\ref{sec:bogomolov}, we prove that
Moravec's conjecture holds for finite $p$-groups
of nilpotency class at most $2p-1$.  See Table~\ref{tab:comparison-bogomolov-bounds} for a detailed comparison. 
Theorem~\ref{thm:exterior-tensor-bounds}
is based on the following general result.

\begin{theorem}\label{thm:higher-central-extension}
Let $p$ be an odd prime and let $G$ be a finite nontrivial  $p$-group.
\begin{enumerate}
    \item If $\cl(G)\le2p-1$, then \(
   \expn G'\mid\expn(G/Z(G)).
\)

\item More generally, for an extension
\(
   1\rightarrow N\rightarrow G\rightarrow Q\rightarrow1
\)
of finite $p$-groups with $N\le Z_s(G)$,
where $s\ge1$, we have
\[
   \expn\g_{s+1}(G)
   \mid
   (\expn Q)^{1+\left\lceil
      \log_p\dfrac{c+1}{2p-1}
   \right\rceil}
\]
with $c=\cl(Q)$.

\end{enumerate}   
\end{theorem}
Our results refine the exponent bounds for
$\g_{s+1}(G)$ previously established by Ellis
\cite{Ellis2001}. See Table~\ref{tab:exponent-comparison}.
Moreover, Vaughan-Lee
\cite[Section~5]{VaughanLee2021counterexamples}
records a conjecture of George Havas that,
for every prime $p>3$, there exists a
four-generator $p$-group of exponent $p$
and nilpotency class $2p-1$ whose Schur
multiplier has exponent greater than $p$.

\begin{theorem}\label{thm:critical-three-generators}
Let $p$ be an odd prime and let $G$ be a finite $p$-group with $\cl(G)\le2p-1$.  
\begin{enumerate}
    \item If $d(G)\le3$, then \(
   \expn\M(G)\mid\expn G.
\)
\item If $\expn G=p$ and  $d(G)=4,$ then \[
   \M(G)\cong C_{p^2}^{\,b}\times C_p^{\,a}
\]
for some $a\ge0$ and $b\in\{0,1\}$.
\end{enumerate}
\end{theorem}

Thus Havas' conjecture can be reformulated as the following assertion about the structure of the Schur multiplier.
\begin{conjecture}
\label{conj:critical-sharpness}
For each prime $p>3$, let $G_p$ be the largest quotient of
\[
    \langle a,b,c,d\mid [b,a]=[d,c]\rangle
\]
of exponent $p$ and nilpotency class at most $2p-1$. Then
\[
    \M(G_p)\cong C_{p^2}\times C_p^{\,a_p}
\]
for some nonnegative integer $a_p$.
\end{conjecture}

The paper is organized as follows.  Section~\ref{sec:preliminaries} collects the preliminary results, including the Hall collection formula and the lower central series.
Section~\ref{sec:induction} develops the two-generator reduction method used in Theorem \ref{thm:higher-central-extension}  $(1)$  and the critical-class analysis in Section~\ref{sec:critical-class} .  Section~\ref{sec:central} proves
Theorem \ref{thm:higher-central-extension} $(1)$ and Theorem \ref{thm:exterior-tensor-bounds} $(1)$.
Sections~\ref{sec:higher-central} and~\ref{sec:schur}
prove Theorems~\ref{thm:higher-central-extension}(2)
and~\ref{thm:exterior-tensor-bounds}(2), respectively.
Section~\ref{sec:critical-class} proves Theorem  \ref{thm:critical-three-generators}.

\FloatBarrier

\section{Notation and preliminary results}\label{sec:preliminaries}

\subsection{Commutators, lower central series, and power subgroups}

We use the conventions
\[
   [x,y]=x^{-1}y^{-1}xy,\qquad x^y=y^{-1}xy,
\]
and, for $r\ge3$, define recursively \(
   [x_1,\ldots,x_r]=[[x_1,\ldots,x_{r-1}],x_r].
\) Thus \(
   [xy,z]=[x,z]^y[y,z]\) and \([x,yz]=[x,z][x,y]^z.
\) 

The lower central series is defined by \(
   \g_1(G)=G,\g_{i+1}(G)=[\g_i(G),G].
\) Thus $G'=\g_2(G)$ and $[\g_i(G),{}_sG]=\g_{i+s}(G)$. Denote $\cl(1)=0.$ And for a nontrivial nilpotent group $G$, $\cl(G)=c$ means $\g_{c+1}(G)=1$ and $\g_c(G)\ne1$. For a group $G$ and a positive integer $q$, write $G^q=\langle g^q:g\in G\rangle$.

\begin{lemma}[{\cite[Theorem~5.1.11(i)]{Robinson1996}}]\label{lem:lower-central-weight}
For every group $G$ and all $r,s\ge1$,
\[
   [\g_r(G),\g_s(G)]\le\g_{r+s}(G).
\]
\end{lemma} It follows that
\begin{equation}\label{eq:derived-weight-consequences}
   \g_i(G')\le\g_{2i}(G),\qquad
   \g_i(\g_3(G))\le\g_{3i}(G).
\end{equation}
More generally, induction gives $\g_i(\g_j(G))\le\g_{ij}(G)$. If $u_j\in\g_{r_j}(G)$ for $1\le j\le m$, then repeated application of Lemma~\ref{lem:lower-central-weight} also gives \(
   [u_1,\ldots,u_m]\in\g_{r_1+\cdots+r_m}(G).
\)

A commutator word of \emph{weight} $d$ is an expression
$[x_1,\ldots,x_d]$, where $d\ge1$,
$[x_1]=x_1$, and the $x_i$ are variables.
If the distinct variables occurring in the word
are $y_1,\ldots,y_m$, its \emph{multiweight}
is $(n_1,\ldots,n_m)$, where $n_i$ is the number
of occurrences of $y_i$.
Thus $d=n_1+\cdots+n_m$.


\begin{lemma}[{\cite[Chapter~5, Proposition~5.3(5),(6)]{MagnusKarrassSolitar1976}}]
\label{lem:multilinearity}
Let $G$ be a group, and let $c(x_1,\ldots,x_d)$ be a commutator word in which each of the $d\ge1$ formal variables occurs once. Let $1\le i\le d$. Then, for $x,y\in G$,
\[
c(\ldots,\underbrace{xy}_{i\text{-th position}},
\ldots) \equiv
c(\ldots,\underbrace{x}_{i\text{-th position}},\ldots)\cdot
c(\ldots,\underbrace{y}_{i\text{-th position}},\ldots)
\pmod{\g_{d+1}(G)}.
\]
\end{lemma}
Thus, for every integer $m$ and every $i$,
\[
   c(x_1,\ldots,x_i^m,\ldots,x_d)
   \equiv c(x_1,\ldots,x_i,\ldots,x_d)^m
   \pmod{\g_{d+1}(G)}.
\]

\subsection{Hall collection and regular
\texorpdfstring{$p$}{p}-groups}

We recall the Hall collection formula and
some properties of regular $p$-groups.

\begin{proposition}[{\cite[Theorem~12.3.1]{Hall1959}}]
\label{prop:hall-collection}
Let $P$ be a finite $p$-group and let
$g_1,\ldots,g_m\in P$. Then:
\begin{enumerate}
\item For every integer $n\ge1$,
\[
   (g_1\cdots g_m)^n
   =
   g_1^n\cdots g_m^n
   \prod_j h_j^{e_j(n)},
\]

where the $h_j$ are basic commutators in
$g_1,\ldots,g_m$ of weight at least $2$,
taken in Hall order.
If $h_j$ has weight $r$, then $e_j$ is an
integral linear combination of
$\binom n1,\ldots,\binom nr$.

\item In particular, if $n=p$ and $h_j$ has
weight $r$ with $2\le r<p$, then
$p\mid e_j(p)$.
\end{enumerate}
\end{proposition}

\begin{proposition}[P.~Hall
{\cite[Corollary~4.13 and Theorem~4.26]{Hall1934}}]
\label{prop:regular-pgroups}
Let $P$ be a finite $p$-group with
$\cl(P)<p$, and let $q=p^n$ with $n\ge1$.
Then:
\begin{enumerate}
\item If $\expn P'\mid q$, then
$(uv)^q=u^qv^q$ for all $u,v\in P$.

\item If $P$ is generated by elements whose
orders divide $q$, then $\expn P\mid q$.
\end{enumerate}
\end{proposition}

\begin{proof}
For (1), Proposition~\ref{prop:hall-collection}
gives the result, since
$q\mid\binom{q}{k}$ for $1\le k<p$
and $\expn P'\mid q$.
For (2), recall that $P$ is regular by
\cite[Corollary~4.13]{Hall1934}.
The conclusion follows from
\cite[Theorem~4.26]{Hall1934}.
\end{proof}

\begin{lemma}\label{lem:power-cosets}
Let $H$ be a finite $p$-group with $\cl(H)\le2p$, let $q=p^n$ with $n\ge1$, and suppose that $\expn\g_3(H)\mid q$. Then \(
   (uh)^q=u^q
\) for every $u\in H'$ and $h\in\g_3(H)$.
\end{lemma}

\begin{proof}
Let $S=\langle u,h\rangle$. Since $[u,h]\in\g_5(H)\trianglelefteq H$, the image of $S$ in $H/\g_5(H)$ is abelian. Hence $S'\le\g_5(H)$. Moreover, $S\le\g_2(H)$, and induction using Lemma~\ref{lem:lower-central-weight} gives $\g_i(S)\le\g_{2i+1}(H)$ for $i\ge2$. Thus $\g_p(S)\le\g_{2p+1}(H)=1$, so $\cl(S)<p$. Since $S'\le\g_5(H)\le\g_3(H)$, we have $\expn S'\mid q$. Also $h^q=1$. Proposition~\ref{prop:regular-pgroups}(1), applied to $S$, gives \(
   (uh)^q=u^qh^q=u^q.\)
\end{proof}

\subsection{Powers in commutator words}

For $V=\langle a,b\rangle$, write $A_d(V)=\g_d(V)/\g_{d+1}(V)$, with $A_d(V)$ written additively.  
Then
\[
   A_1(V)=\langle aV',bV'\rangle,
\]
and for $d\ge2$,
\[
   A_d(V)=\big\langle [x_1,\ldots,x_d]\g_{d+1}(V):
      x_1,\ldots,x_d\in\{a,b\}\big\rangle.
\]
\begin{lemma}\label{lem:scalar-factor}
Let $V=\langle a,b\rangle$, and let $\sigma_{r,s}$ be an endomorphism of $V$ satisfying $\sigma_{r,s}(a)=a^r$ and $\sigma_{r,s}(b)=b^s$ for some integers $r,s$. If $w$ is a commutator word of multiweight $(i,j)$ in $a,b$, then
\[
   \sigma_{r,s}(w)
   \equiv w^{r^is^j}
   \pmod{\g_{i+j+1}(V)}.
\]
\end{lemma}

\begin{proof}
Write $
   w=c(x_1,\ldots,x_{i+j}),
$ where each $x_k$ is either $a$ or $b$, with exactly $i$ of the $x_k$ equal to $a$ and exactly $j$ equal to $b$. Then
\[
   \sigma_{r,s}(w)
   =
   c(x_1^{m_1},\ldots,x_{i+j}^{m_{i+j}}),
\]
where $m_k=r$ if $x_k=a$, or $m_k=s$ if $x_k=b$. By Lemma~\ref{lem:multilinearity}, for every integer $m$ and every position $k$,
\[
   c(x_1,\ldots,x_k^m,\ldots,x_{i+j})
   \equiv
   c(x_1,\ldots,x_k,\ldots,x_{i+j})^m
   \pmod{\g_{i+j+1}(V)}.
\]
Applying this congruence successively to the $i+j$ occurrences gives
\[
\begin{aligned}
   \sigma_{r,s}(w)
   =
   c(x_1^{m_1},\ldots,x_{i+j}^{m_{i+j}})
   \equiv
   c(x_1,\ldots,x_{i+j})^{m_1\cdots m_{i+j}}
   \pmod{\g_{i+j+1}(V)}.
\end{aligned}
\]
Since $i$ of the integers $m_k$ are equal to $r$ and the remaining $j$ are equal to $s$, we have
$
   m_1\cdots m_{i+j}=r^is^j.
$ This completes the proof.
\end{proof}

\begin{remark}
When $r=s$, Lemma~\ref{lem:scalar-factor} shows that $\sigma_{r,r}$ induces multiplication by $r^d$ on $\g_d(V)/\g_{d+1}(V)$ for every $d\ge1$. This special case is used in the proof of Theorem~\ref{thm:higher-central-extension}(1), while the two-parameter form is used in Section~\ref{sec:critical-class}.
\end{remark}

\subsection{Schur covers, exterior squares, and CP covers}

Let $G=F/R$ be a free presentation, and write $\M(G)=H_2(G,\Z)$.
A central extension $1\to N\to E\to G\to1$ is \emph{stem} if
$N\le E'$.  A \emph{Schur cover} is a stem central extension whose
kernel is isomorphic to $\M(G)$.

\begin{proposition}\label{prop:schur-cover-facts}
Let $G=F/R$ be a finite group.
\begin{enumerate}
\item[(1)] $\M(G)\cong(F'\cap R)/[F,R]$
\cite[Theorem~2.4.6]{Karpilovsky1987}.
\item[(2)] There exists a Schur cover
$1\to\M(G)\to\widetilde G\to G\to1$
\cite[Theorem~2.1.4]{Karpilovsky1987}.
\item[(3)] If $G$ is a finite $p$-group, then
$\M(G)$ and every Schur cover $\widetilde G$
are finite $p$-groups
\cite[Theorems~2.1.2 and~2.1.4]{Karpilovsky1987}.
\item[(4)] For every Schur cover $\widetilde G$ of $G$,
we have $\widetilde G'\cong G\wed G$
\cite[Corollary~2]{BrownJohnsonRobertson1987}.
\end{enumerate}
\end{proposition}

We use the nonabelian tensor square $G\otimes G$ introduced by
Brown--Loday \cite{BrownLoday1987}, with diagonal subgroup
$\nabla(G)=\langle g\otimes g:g\in G\rangle$, and define
$G\wed G=(G\otimes G)/\nabla(G)$.  
For a free presentation $G=F/R$,
\cite[Theorem~5 and its proof]{Ellis1987Exterior}
gives $G\wed G\cong F'/[F,R]$
and the exact sequence
\begin{equation}\label{eq:exterior-exact}
   1\longrightarrow\M(G)\longrightarrow G\wed G
   \longrightarrow G'\longrightarrow1.
\end{equation}
For any group $E$, let $\Kset(E)=\{[x,y]:x,y\in E\}$.
Define $N_0=\langle\Kset(F)\cap R\rangle^F$, where the superscript
denotes normal closure.  Write $\B(G)$ for the  Bogomolov
multiplier.  The multiplier $\B(G)$ and curly exterior square have the presentation forms
\begin{equation}\label{eq:curly-models}
   \B(G)\cong(F'\cap R)/N_0,\qquad G\curlyw G\cong F'/N_0.
\end{equation}
By \cite[Proposition~3.8 and the proof of
Theorem~4.2]{Moravec2012Unramified},
we have the exact sequence
\[
   1\longrightarrow\B(G)\longrightarrow G\curlyw G
   \longrightarrow G'\longrightarrow1.
\]
Since $[F,R]\le N_0\le F'\cap R$, the group $\B(G)$ is a quotient
of $\M(G)$.

A central extension is \emph{commutativity preserving} (CP) if each
commuting pair in the quotient has commuting lifts.  A \emph{CP cover}
is a stem central CP extension with kernel isomorphic to $\B(G)$.

\begin{proposition}\label{prop:cp-package}
Let
$1\to N\to E\overset{\pi}{\to}G\to1$ be a central extension.
Then the following statements hold.
\begin{enumerate}
\item[(1)] The extension is CP if and only if
$N\cap\Kset(E)=\{1\}$
\cite[Proposition~3.3]{JezernikMoravec2018}.
\item[(2)] If $G$ is finite, then $G$ admits a CP cover
$1\to\B(G)\to G^*\to G\to1$
\cite[Theorem~4.2(1),(2)]{JezernikMoravec2018}.
\item[(3)] If the extension $E\overset{\pi}{\to}G$ is CP, then
$\pi$ induces an isomorphism $E/Z(E)\cong G/Z(G)$.
\item[(4)] If $G$ is finite and $G^*\to G$ is a CP cover,
then $(G^*)'\cong G\curlyw G$
\cite{Moravec2019Bogomolov}.
If $G$ is a finite $p$-group, then $G^*$ is
a finite $p$-group.
\end{enumerate}
\end{proposition}

\begin{proof}

For~(3), if $\pi(e)\in Z(G)$, then
$[e,f]\in N\cap\Kset(E)$ for every $f\in E$.
By part~(1), we have $e\in Z(E)$, and hence
$\pi^{-1}(Z(G))=Z(E)$.
For the last assertion in~(4),
$\B(G)$ is a quotient of $\M(G)$, so
Proposition~\ref{prop:schur-cover-facts}(3)
and the exact sequence of the CP cover
show that $G^*$ is a finite $p$-group.
\end{proof}

\section{Two-generator reductions}\label{sec:induction}

Fix an odd prime $p$. To prove that $\expn G'\mid\expn(G/Z(G))$ for finite $p$-groups
of class at most $2p-1$, we reduce the problem to a commutator calculation in a two-generator group, using a construction similar to that of Vaughan-Lee \cite[Section~5]{VaughanLee2021II}.

\subsection{A two-generator construction}

We give the two-generator construction firstly.  Let $q=p^n$ with $n\ge1$, let $k\ge2$, and let $F=\langle A,B\rangle$ be free of rank two.  Define
\begin{equation}\label{eq:two-generator-U}
   R=F^q\g_k(F),\qquad U=F/[F,R].
\end{equation}
Let $a=A[F,R]$ and $b=B[F,R]$.  Since $R/[F,R]\le Z(U)$, define
\begin{equation}\label{eq:two-generator-V}
   C=\langle a^q,b^q\rangle\le Z(U),\qquad V=U/C.
\end{equation}
Write $\bar a=aC$ and $\bar b=bC$.

\begin{proposition}\label{prop:universal-reduction}
The groups $U$ and $V$ satisfy the following properties.
\begin{enumerate}
\item[(1)] If $H$ is a group satisfying $H^q\le Z(H)$
and $\cl(H)\le k$, then for any $x,y\in H$
there exists a surjective homomorphism
\[
   \alpha:U\to\langle x,y\rangle
\]
such that $\alpha(a)=x$ and $\alpha(b)=y$.

\item[(2)] The quotient map $\pi:U\to V$ restricts to an isomorphism
\[
   \pi|_{U'}:U'\to V'.
\]
\item[(3)] $V$ is a finite $p$-group satisfying
\[
   \bar a^q=\bar b^q=1,\qquad
   \cl(V)\le k,\qquad
   V^q\le Z(V).
\]

\item[(4)] For all integers $r,s$ with $p\nmid rs$,
there is an automorphism $\sigma_{r,s}\in\Aut(V)$
such that
\[
   \sigma_{r,s}(\bar a)=\bar a^r
   \qquad\text{and}\qquad
   \sigma_{r,s}(\bar b)=\bar b^s.
\]
\end{enumerate}
\end{proposition}

\begin{proof}
For~$(1)$, let $\beta:F\to\langle x,y\rangle$ be defined by $\beta(A)=x$ and $\beta(B)=y$.  Since $H^q\le Z(H)$ and $\cl(H)\le k$, we have $\beta(F^q)\le Z(H)$ and $\beta(\g_k(F))\le\g_k(H)\le Z(H)$.  Hence $\beta(R)\le Z(H)$, so $[F,R]\le\ker\beta$.  Thus $\beta$ factors through $U$ to a surjective homomorphism $\alpha:U\to\langle x,y\rangle$ satisfying $\alpha(a)=x$ and $\alpha(b)=y$.

For~(2), since $[F,R]\le F'$, we have $U'=F'/[F,R]$, and hence
\[
   U/U'\cong
   \frac{F/[F,R]}{F'/[F,R]}
   \cong F/F'\cong\Z^2.
\]
Under this identification, $aU'$ and $bU'$ correspond to $(1,0)$ and $(0,1)$.  Since $a^q$ and $b^q$ are central, every element of $C$ has the form $a^{qi}b^{qj}$ for some $i,j\in\Z$, and $a^{qi}b^{qj}U'=(qi,qj)$.  Hence $a^{qi}b^{qj}\in U'$ implies $i=j=0$, and $C\cap U'=1$. Since $\pi$ is surjective, $\pi(U')=V'$, while $\ker(\pi|_{U'})=U'\cap C=1$.  Hence $\pi|_{U'}:U'\to V'$ is an isomorphism.

For~(3), the relations $\bar a^q=\bar b^q=1$ hold in $V$ by construction, while
\[
   \g_{k+1}(F)=[\g_k(F),F]\le[R,F]=[F,R],
\]
so $\cl(V)\le k$. The group $V/V'$ is generated by $\bar aV'$ and $\bar bV'$, both of order dividing $q$, and hence is a finite $p$-group.  For $i\ge2$, the abelian group $\g_i(V)/\g_{i+1}(V)$ is generated by the elements  $ t\g_{i+1}(V)$ where $t=[x_1,\ldots,x_i]$  with $x_1,\ldots,x_i\in\{\bar a,\bar b\}.$ For $1\le j\le i$, Lemma~\ref{lem:multilinearity} gives
\[
   1=[x_1,\ldots,x_j^q,\ldots,x_i]
   \equiv t^q\pmod{\g_{i+1}(V)}.
\]
Thus $\g_i(V)/\g_{i+1}(V)$ is a finite abelian $p$-group of exponent dividing $q$.  Since $\g_{k+1}(V)=1$, it follows that $V$ is a finite $p$-group. If $u=f[F,R]\in U$, then
\[
   u^q=f^q[F,R]\in R/[F,R]\le Z(U),
\]
since $F^q\le R$.  Hence $U^q\le Z(U)$, and therefore $V^q\le Z(V)$.

For~(4), let $\varphi_{r,s}$ be the endomorphism of $F$ defined by $\varphi_{r,s}(A)=A^r$ and $\varphi_{r,s}(B)=B^s$.  The subgroups
$F^q$ and $\g_k(F)$ are fully invariant, so $\varphi_{r,s}$ preserves $R$ and hence $[F,R]$.  Therefore $\varphi_{r,s}$ induces an
endomorphism $\overline\varphi_{r,s}$ of $U$. Since
\[
   \overline\varphi_{r,s}(a^q)=a^{rq}
   \qquad\text{and}\qquad
   \overline\varphi_{r,s}(b^q)=b^{sq},
\]
we have $\overline\varphi_{r,s}(C)\le C$, and hence $\overline\varphi_{r,s}$ induces an endomorphism $\sigma_{r,s}$ of $V$ satisfying $\sigma_{r,s}(\bar a)=\bar a^r$ and $\sigma_{r,s}(\bar b)=\bar b^s$. Since $r$ and $s$ are coprime to $q$, we can choose $r',s'\in\Z$ such that $rr'\equiv1\pmod q$ and $ss'\equiv1\pmod q$.  As $\bar a^q=\bar b^q=1$, the maps $\sigma_{r,s}$ and $\sigma_{r',s'}$ are inverse on both generators of $V$.  Since $V=\langle\bar a,\bar b\rangle$, they are inverse automorphisms.
\end{proof}

\subsection{The third term of the lower central series}

For $c\ge1$, let $\calP(c)$ denote the assertion that $\expn G'\mid\expn(G/Z(G))$ for every finite
$p$-group $G$ with $\cl(G)\le c$. To prove Theorem~\ref{thm:higher-central-extension}(1), we use the result of Antony, Komma and Thomas as the base case $\calP(p+1)$ and establish $\calP(c)\Rightarrow\calP(c+1)$ for $p+1\le c<2p-1$.

\begin{lemma}\label{lem:gamma3}
Let $c$ be an integer with $p+1\le c\le2p-1$,
and suppose that $\calP(c)$ holds.
Let $q$ be a power of $p$, and let $V$ be a finite
$p$-group with $\cl(V)\le c+1$ and $V^q\le Z(V)$.
Then $\expn\g_3(V)\mid q$.
\end{lemma}

\begin{proof}
Let $u\in V'$ and $v\in V$, and let $L=\langle u,v\rangle$. Since $[u,v]\in\g_3(V)\trianglelefteq V$, we have $L'\le\g_3(V)$. It follows by induction that $\g_i(L)\le\g_{i+1}(V)$ for all $i\ge2$. As $\cl(V)\le c+1$, we have $\g_{c+1}(L)\le\g_{c+2}(V)=1$, so $\cl(L)\le c$. Since $V^q\le Z(V)$, we have $l^q\in L\cap Z(V)\le Z(L)$ for every $l\in L$. Thus $\expn(L/Z(L))\mid q$. By $\calP(c)$, we obtain $\expn L'\mid q$, and hence $[u,v]^q=1$.

The commutators $[u,v]$, with $u\in V'$ and $v\in V$, generate $\g_3(V)$. By \eqref{eq:derived-weight-consequences}, \(
\g_p(\g_3(V))\le\g_{3p}(V)=1,
\) since $\cl(V)\le c+1\le2p$. Thus $\g_3(V)$ is a finite $p$-group of class less than $p$, generated by elements of order dividing $q$. Proposition~\ref{prop:regular-pgroups}(2) gives $\expn\g_3(V)\mid q$.
\end{proof}

\subsection{Commutator powers}
\begin{lemma}[Antony--Komma--Thomas
{\cite[Theorem~2.4]{AntonyKommaThomas2022}}]
\label{lem:akt-base}
Let $p$ be a prime and let $P$ be a finite $p$-group. If $\cl(P)\le p+1$, then
\[
   \expn P'\mid\expn(P/Z(P)).
\]
\end{lemma}

\begin{proposition}\label{prop:power-scaling}
Let $V$ be the group constructed in
\eqref{eq:two-generator-U}--\eqref{eq:two-generator-V}
with $k\le2p$, and suppose that
$\expn\g_3(V)\mid q$.
Set $w=[\bar b,\bar a]$ and $z=w^q$.
Then $z\in\g_{p+2}(V)$, and
$
   \sigma_{r,s}(z)=z^{rs}
$
for all integers $r,s$ with $p\nmid rs$.
\end{proposition}

\begin{proof}
Let $\rho:V\to\overline V=V/\g_{p+2}(V)$ be the quotient map. Then $\cl(\overline V)\le p+1$. By Proposition~\ref{prop:universal-reduction}(3), $V^q\le Z(V)$, and hence $\overline V^q\le Z(\overline V)$. Thus \(
   \expn(\overline V/Z(\overline V))\mid q.
\) By Lemma~\ref{lem:akt-base}, $\expn\overline V'\mid q$.  Since $w\in V'$,  $\rho(w)^q=1$.  Therefore \(
   z=w^q\in\ker\rho=\g_{p+2}(V).
\) Now let $p\nmid rs$.  Lemma~\ref{lem:scalar-factor} gives
\[
   [\bar b^s,\bar a^r]
   \equiv
   [\bar b,\bar a]^{rs}
   =
   w^{rs}
   \pmod{\g_3(V)}.
\]
Thus $\sigma_{r,s}(w)=w^{rs}h$ for some $h\in\g_3(V)$.  Lemma~\ref{lem:power-cosets} gives
\[
   \sigma_{r,s}(z)
   =\sigma_{r,s}(w)^q
   =(w^{rs}h)^q
   =w^{rsq}
   =z^{rs}.
\]
\end{proof}

\begin{lemma}\label{lem:power-depth}
Under the hypotheses and notation of
Proposition~\ref{prop:power-scaling},
we have $z\in\g_{2p}(V)$. In particular, if $k\le2p-1$, then $z=1$.
\end{lemma}

\begin{proof}
Choose a positive integer $r$ whose residue modulo $p$ generates $\F_p^\times$, and write $\sigma=\sigma_{r,r}$. By Proposition~\ref{prop:power-scaling}, \(z\in\g_{p+2}(V)\) and  \(
   \sigma(z)=z^{r^2}.
\) We claim that
\[
   z\in\g_d(V)
   \quad\Longrightarrow\quad
   z\in\g_{d+1}(V)
\]
for every $d$ with $p+2\le d\le2p-1$. Fix such a $d$ and suppose that $z\in\g_d(V)$. 
Let $\overline z=z\g_{d+1}(V)\in A_d(V)$.
The group $A_d(V)$ is generated by the images of weight-$d$ commutators in the distinguished generators $\bar a,\bar b$. By Lemma~\ref{lem:scalar-factor}, the automorphism $\sigma=\sigma_{r,r}$ therefore acts on $A_d(V)$ as multiplication by $r^d$. Hence \(
   \sigma(\overline z)=r^d\overline z.
\) On the other hand, $\sigma(z)=z^{r^2}$ gives \(
   \sigma(\overline z)=r^2\overline z.
\) Consequently
\[
   (r^d-r^2)\overline z=0.
\]
Since the residue of $r$ modulo $p$ has multiplicative order $p-1$, we have
\[
   p\mid(r^d-r^2)
   \quad\Longleftrightarrow\quad
   r^{d-2}\equiv1\pmod p
   \quad\Longleftrightarrow\quad
   d\equiv2\pmod{p-1}.
\]
Thus $p-1$ does not divide $d-2$ by the range of $d$, and hence $p\nmid(r^d-r^2)$. 
Since $V$ is a finite $p$-group, $\overline z$
has $p$-power order. As $p\nmid(r^d-r^2)$,
we obtain $\overline z=0$, and hence
$z\in\g_{d+1}(V)$.

Applying the preceding claim inductively for
$d=p+2,\ldots,2p-1$, starting with
$z\in\g_{p+2}(V)$, we obtain
$z\in\g_{2p}(V)$.
Finally, if $k\le2p-1$, then
Proposition~\ref{prop:universal-reduction}(3)
gives $\cl(V)\le k\le2p-1$, so
$\g_{2p}(V)=1$. Hence $z=1$.
\end{proof}

\section{The main central-factor theorem}\label{sec:central}

\begin{proof}[Proof of Theorem~\ref{thm:higher-central-extension} $(1)$]
We prove $\calP(c)$ for $p+1\le c\le2p-1$ by induction on $c$. The base case $\calP(p+1)$ is Lemma~\ref{lem:akt-base}. Suppose that $p+1\le c<2p-1$ and that $\calP(c)$ holds. Let $X$ be a finite $p$-group with $\cl(X)\le c+1$, and let $q=\expn(X/Z(X))$. If $q=1$, then $X=Z(X)$ and $X'=1$. We may therefore assume that $q=p^n$ with $n\ge1$. Fix $x,y\in X$. Since $X^q\le Z(X)$ and $\cl(X)\le c+1$, Proposition~\ref{prop:universal-reduction} applies with $k=c+1$. Let
\[
   \alpha:U\longrightarrow\langle x,y\rangle
   \qquad\text{and}\qquad
   \pi:U\longrightarrow V
\]
be the maps in that proposition, and write \(
   z=[\bar a,\bar b]^q.
\)
By Proposition~\ref{prop:universal-reduction}(3), we have $\cl(V)\le c+1$ and $V^q\le Z(V)$. Thus Lemma~\ref{lem:gamma3} gives $\expn\g_3(V)\mid q$. Hence Lemma~\ref{lem:power-depth} gives $z=1$. By Proposition~\ref{prop:universal-reduction}(2), the restriction $\pi|_{U'}:U'\to V'$ is injective. Hence $[b,a]^q=1$ in $U$. Applying $\alpha$, which sends $a,b$ to $x,y$, gives $[y,x]^q=1$, and hence $[x,y]^q=1$.

Finally, \eqref{eq:derived-weight-consequences} gives \(
   \g_p(X')\le\g_{2p}(X)=1,
\) because $\cl(X)\le c+1\le2p-1$. Moreover, $X'$ is generated by the commutators $[x,y]$, all of which have order dividing $q$. By Proposition~\ref{prop:regular-pgroups}(2),
\[
   \expn X'\mid \expn(X/Z(X))=q.
\]
Therefore $\calP(c+1)$ holds, completing the induction.
\end{proof}

\begin{proof}[Proof of Theorem \ref{thm:exterior-tensor-bounds}$(1)$]
Choose a Schur cover
\[
   1\longrightarrow\M(G)\longrightarrow\widetilde G
   \overset{\pi}{\longrightarrow}G\longrightarrow1,
   \qquad \M(G)\le Z(\widetilde G)\cap\widetilde G'.
\]
By Proposition~\ref{prop:schur-cover-facts}(3), $\widetilde G$ is a finite $p$-group.  Since $\g_{2p-1}(G)=1,$ we have \(
   \g_{2p-1}(\widetilde G)\le\M(G)\le Z(\widetilde G).
\) Thus $\cl(\widetilde G)\le2p-1$. So Theorem~\ref{thm:higher-central-extension} gives
\[
   \expn\M(G)
   \mid\expn\widetilde G'
   \mid\expn(\widetilde G/Z(\widetilde G))
   \mid\expn G.
\]

\end{proof}

\begin{remark}\label{rem:degree-2p}
For the chosen $r$, the first integer $d>p+1$
satisfying $r^d\equiv r^2\pmod p$ is $d=2p$.
Thus the induction in Lemma~\ref{lem:power-depth}
gives $z\in\g_{2p}(V)$ but does not imply
$z=1$ when $\cl(V)=2p$.
In Section~\ref{sec:critical-class}, we use
the two-parameter automorphisms from
Proposition~\ref{prop:universal-reduction}(4)
to analyse the commutators of weight $2p$.
\end{remark}

\section{Higher central-factor bounds}\label{sec:higher-central}

In this section, we extend the divisibility $\expn G'\mid\expn(G/Z(G))$ from Theorem~\ref{thm:higher-central-extension}(1) to $\expn\g_{s+1}(G)\mid\expn(G/Z_s(G))$ Theorem~\ref{thm:higher-central-extension}(2). For a group $G$, write
\[
   1=Z_0(G)\le Z_1(G)\le Z_2(G)\le\cdots
\]
for the upper central series. Thus $Z_1(G)=Z(G)$ and
\[
   Z_{i+1}(G)/Z_i(G)=Z(G/Z_i(G)).
\]

\subsection{Higher-central exponent divisibility}

\begin{lemma}\label{lem:descending-exponent-step}
Let $p$ be an odd prime, $G$ a finite $p$-group, $k$ a positive integer, and $q$ a power of $p$. If $\cl(G)\le pk$ and $\g_k(G)^q\le Z(G)$, then $
   \expn\g_{k+1}(G)\mid q.
$
\end{lemma}

\begin{proof}

If $q=1$, then $\g_k(G)\le Z(G)$ and hence
$\g_{k+1}(G)=1$.
We may therefore assume $q>1$ and proceed
by descending induction on $k$. For $k\ge\cl(G)$, we have $\g_{k+1}(G)=1$, so the assertion is immediate. Now assume that the assertion holds for $k+1$, and suppose that $\cl(G)\le pk$ and $\g_k(G)^q\le Z(G)$. Since
\[
   \cl(G)\le p(k+1) \quad\text{ and } \quad
   \g_{k+1}(G)^q\le\g_k(G)^q\le Z(G),
\]
the induction hypothesis gives $\expn\g_{k+2}(G)\mid q$. We now use this to prove that $\expn\g_{k+1}(G)\mid q$. Let $a\in\g_k(G)$, let $b\in G$, and let $
   d=[a,b],
   P=\langle a,d\rangle.
$ Since $
   d\in\g_{k+1}(G),
$ we have
\[
   [d,a]
   \in[\g_{k+1}(G),\g_k(G)]
   \le\g_{2k+1}(G).
\]
Since $P=\langle a,d\rangle$, we obtain \(
   P'\le\g_{2k+1}(G)\le\g_{k+2}(G).
\) Thus $\expn P'\mid q$. Moreover, induction on $i$ gives $\g_i(P)\le\g_{ik+1}(G)$ for $i\ge2$. Since $\cl(G)\le pk$, we have $\g_p(P)\le\g_{pk+1}(G)=1$. Thus $\cl(P)<p$, and Proposition~\ref{prop:regular-pgroups}(1) gives $(ad)^q=a^qd^q$. Since $a^q\in\g_k(G)^q\le Z(G)$ and $ad=a^b$,
\[
   a^qd^q=(ad)^q=(a^b)^q=(a^q)^b=a^q,
\]
and hence $d^q=1$. 
The commutators $d=[a,b]$ generate $\g_{k+1}(G)$.
Furthermore, $\g_p(\g_{k+1}(G))
\le\g_{p(k+1)}(G)=1$. Thus $\g_{k+1}(G)$ has class less than $p$, and Proposition~\ref{prop:regular-pgroups}(2) yields $\expn\g_{k+1}(G)\mid q$. This completes the descending induction.
\end{proof}

\begin{theorem}\label{thm:higher-central}
Let $p$ be an odd prime, $G$ a finite $p$-group, and $s$ a positive integer. If $\cl(G)\le 2p+s-2$, then $
   \expn\g_{s+1}(G)\mid\expn(G/Z_s(G)).
$
\end{theorem}

\begin{proof}[Proof of Theorem~\ref{thm:higher-central}]
We induct on $s$.  The case $s=1$ is Theorem~\ref{thm:higher-central-extension}(1). Assume $s\ge2$ and that the assertion holds with $s-1$ in place of $s$. Let $
q=\expn(G/Z_s(G))=p^n.
$ If $q=1$, then $G=Z_s(G)$ and $\g_{s+1}(G)=1$.  We therefore assume $n\ge1$.
We first show that $\g_s(G)^q\le Z(G)$.

Let $\overline G=G/Z(G)$. Since
$
   Z_{s-1}(\overline G)=Z_s(G)/Z(G),
   \cl(\overline G)\le2p+s-3,
$
we have $\overline G/Z_{s-1}(\overline G)\cong G/Z_s(G)$, and the inductive hypothesis gives \(   \expn\g_s(\overline G)\mid \expn(\overline G/Z_{s-1}(\overline G))=\expn(G/Z_s(G))=q.\) Since $\g_s(\overline G)=\g_s(G)Z(G)/Z(G)$, thus $a^q\in Z(G)$ for every $a\in\g_s(G)$. Since \(s\ge 2\), we have
$
   \cl(G)\le 2p+s-2\le ps.
$
Hence Lemma~\ref{lem:descending-exponent-step}, applied with \(k=s\), yields \(
   \expn\g_{s+1}(G)\mid q,
\) as required.
\end{proof}

\subsection{A logarithmic higher-central bound}\label{subsec:higher-central-logarithmic}


Theorem~\ref{thm:higher-central} gives
$\expn\g_{s+1}(G)\mid\expn(G/Z_s(G))$
when $\cl(G/Z_s(G))\le2p-2$.
We now establish a two-parameter bound that
allows us to iterate this estimate and obtain
a logarithmic exponent bound.

\begin{theorem}\label{thm:higher-central-two-parameter}
Let $p$ be an odd prime, $G$ a finite $p$-group,
and $r,t$ positive integers.
If $\cl(G)\le pr+t-1$, then
$
   \expn\g_{r+t}(G)\mid\expn(G/Z_t(G)).
$
\end{theorem}

\begin{proof}
We first consider $t=1$ and write
$q=\expn(G/Z(G))$. If $q=1$, then
$G$ is abelian and the conclusion is immediate.
We may assume $q>1$.
Since \(
   \g_r(G)^q\le G^q\le Z(G)
\) and \(
   \cl(G)\le pr,
\) Lemma~\ref{lem:descending-exponent-step}, applied with \(k=r\), gives \(
   \expn \g_{r+1}(G)\mid q,
\) which proves the case \(t=1\).

Now we induct on $t$. Let $t\ge2$ and assume the assertion holds with $t-1$ in place of $t$. Let $q=\expn(G/Z_t(G))$.  If $q=1$, then $G=Z_t(G)$ and $\g_{r+t}(G)=1$.  Let $\overline G=G/Z(G)$.  Since $Z_{t-1}(\overline G)=Z_t(G)/Z(G)$ and $\cl(\overline G)\le pr+t-2$, the inductive hypothesis gives \(
   \expn\g_{r+t-1}(\overline G)\mid \expn(\overline G/Z_{t-1}(\overline G))= q.
\) This gives
\begin{equation}\label{eq:higher-log-B-central}
   \g_{r+t-1}(G)^q\le Z(G)\cap \g_{r+t-1}(G)\leq Z(\g_{r+t-1}(G)).
\end{equation}
In particular, \(
   \g_{r+t-1}(G)^q\le Z(G).
\) Moreover, since $t\ge2$, \(
   \cl(G)\le pr+t-1\le p(r+t-1).
\) Hence Lemma~\ref{lem:descending-exponent-step}, applied with $k=r+t-1$, gives $
   \expn\g_{r+t}(G)\mid q.
$
\end{proof}

\begin{corollary}\label{cor:higher-central-multiplicative}
Let $p$ be an odd prime, $G$ a finite $p$-group, and $s,r$ positive integers. Set $q=\expn(G/Z_s(G))$. Then
\begin{equation}\label{eq:higher-central-multiplicative}
   \g_{s+r}(G)^q\le\g_{s+pr}(G).
\end{equation}
\end{corollary}

\begin{proof}
Let $\overline G=G/\g_{s+pr}(G)$.  Then $\cl(\overline G)\le s+pr-1$. The natural surjective map $G\twoheadrightarrow\overline G$ maps $Z_s(G)$ into $Z_s(\overline G)$, and hence induces a surjective map \(
   G/Z_s(G)
   \twoheadrightarrow
   \overline G/Z_s(\overline G).
\) Therefore
\[
   \expn(\overline G/Z_s(\overline G))
   \mid
   \expn(G/Z_s(G))
   =q.
\]
Theorem~\ref{thm:higher-central-two-parameter}, applied to $\overline G$ with $t=s$, gives
\[
\begin{aligned}
 \expn (\g_{r+s}(G) /\g_{s+pr}(G))
 = \expn\g_{r+s}(\overline G)
 \mid\expn(\overline G/Z_s(\overline G))\mid q.
\end{aligned}
\]
\end{proof}

For $c=0$, define $ \mu_p(c)=0.$ When $c\geq 1,$ define $\mu_p(c)=  1+\left\lceil\log_p\frac{c+1}{2p-1}\right\rceil. $ Equivalently, for $c\ge1$, the integer $\mu_p(c)$ is the least positive
integer $m$ such that $c+1\le(2p-1)p^{m-1}$.

\begin{theorem}\label{thm:higher-central-logarithmic}
Let $G$ be a finite $p$-group for an odd prime $p$, and let $s$ be a positive integer.
Write $q=\expn(G/Z_s(G))$ and $c=\cl(G/Z_s(G))$. Then $
   \expn\g_{s+1}(G)\mid q^{\mu_p(c)}.
$
\end{theorem}

\begin{proof}
If $c=0$, then $G=Z_s(G)$ and $\g_{s+1}(G)=1$.  Assume $c\ge1$.
Since $\g_{c+1}(G)\le Z_s(G)$, we have
\begin{equation}\label{eq:higher-log-terminal}
   \g_{c+s+1}(G)=1.
\end{equation}
Let $G_0=G/\g_{s+2p-1}(G)$.  Then $\cl(G_0)\le2p+s-2$, and the natural map $G/Z_s(G)\twoheadrightarrow G_0/Z_s(G_0)$ is surjective. Theorem~\ref{thm:higher-central} gives
\begin{equation}\label{eq:higher-log-first-step}
   \g_{s+1}(G)^q\le\g_{s+2p-1}(G).
\end{equation}
Let $m=\mu_p(c)$.  If $m=1$, then $c+1\le2p-1$, and \eqref{eq:higher-log-terminal} and \eqref{eq:higher-log-first-step} give the assertion.  Assume $m\ge2$, and let $r_i=(2p-1)p^{i-1}$ for $1\le i\le m$.  Repeated application of Corollary~\ref{cor:higher-central-multiplicative} gives $\g_{s+r_i}(G)^q\le\g_{s+r_{i+1}}(G)$ for $1\le i<m$. Together with \eqref{eq:higher-log-first-step}, this yields
$
   \g_{s+1}(G)^{q^m}\le\g_{s+r_m}(G).
$
By the definition of $m$, we have $r_m\ge c+1$.  Hence \eqref{eq:higher-log-terminal} gives $\g_{s+r_m}(G)=1$, proving the assertion.
\end{proof}

\begin{proof}[Proof of Theorem~\ref{thm:higher-central-extension}]
The quotient map $Q=G/N\twoheadrightarrow G/Z_s(G)$ is surjective.  Let
\[
   q=\expn(G/Z_s(G)) \quad \text{and } \quad d=\cl(G/Z_s(G)).
\]
Then $q\mid\expn Q$ and $d\le c$.  Theorem~\ref{thm:higher-central-logarithmic} gives
\[
   \expn\g_{s+1}(G)
   \mid q^{\mu_p(d)}
   \mid(\expn Q)^{\mu_p(c)}.
\]
\end{proof}

\begin{remark}\label{rem:ellis-higher-comparison}

We compare Theorem~\ref{thm:higher-central-extension}
with the bound obtained from
Ellis~\cite[Theorem~B(i)]{Ellis2001}. Let
\[
   1\longrightarrow N\longrightarrow G\longrightarrow Q\longrightarrow1
\]
be an extension of finite $p$-groups with $N\le Z_s(G)$, and let \(
   c=\cl(Q).
\) The corresponding exponent bounds are summarized in the following table.
\begin{table}[ht]
\centering
\caption{Exponent bounds for higher central extensions}
\label{tab:exponent-comparison}
\renewcommand{\arraystretch}{1.5}
\begin{tabularx}{\textwidth}{@{}l >{\centering\arraybackslash}X@{}}
\toprule
Source & Bound \\
\midrule

Ellis~\cite[Theorem~B(i)]{Ellis2001}
&
\(
   \expn\g_{s+1}(G)
   \mid
   (\expn Q)^{\lceil c/2\rceil}.
\)
\\

\addlinespace

Theorem~\ref{thm:higher-central-extension}
&
\(
   \expn\g_{s+1}(G)
   \mid
   (\expn Q)^{\,1+\left\lceil
      \log_p\dfrac{c+1}{2p-1}
   \right\rceil}.
\)
\\

\bottomrule
\end{tabularx}
\end{table}

For the special case $\expn(G/Z(G))=p$,
let $Q=G/Z(G)$ and $d=\cl(G)$.
Since $\cl(Q)=d-1$, the results of
Komma--Thomas~\cite[Theorem~4.1]{KommaThomas2025}
and Theorem~\ref{thm:higher-central-extension}
give the following bounds.

\begin{table}[ht]
\centering
\caption{Exponent bound for $\expn(G/Z(G))=p$}
\label{tab:comparison-exp-p}
\renewcommand{\arraystretch}{1.5}
\begin{tabularx}{\textwidth}{@{}l >{\centering\arraybackslash}X@{}}
\toprule
 Source & Bound \\
\midrule

Komma--Thomas~\cite[Theorem~4.1]{KommaThomas2025}
&
\(
    \expn G'\mid p^r \text{ for }  d\le p^{r} 
\)
\\

\addlinespace

Theorem~\ref{thm:higher-central-extension}
&
\(
   \expn G'\mid p^r \text{ for }  d\le (2p-1)p^{r-1} 
\)
\\
\bottomrule
\end{tabularx}
\end{table}

\end{remark}

\section{Multiplier, exterior and tensor consequences}\label{sec:schur}
We now apply the higher-central bounds to
Schur covers and CP covers to bound the exponents
of nonabelian exterior and tensor squares,
Schur multipliers, and Bogomolov multipliers.
\subsection{Exterior and tensor-square bounds}\label{subsec:variants}

\begin{proof}[Proof of Theorem~\ref{thm:exterior-tensor-bounds} (2)]
Choose a Schur cover
\[
   1\longrightarrow\M(G)\longrightarrow\widetilde G
   \longrightarrow G\longrightarrow1.
\]
Since $\M(G)\le Z(\widetilde G)$, Theorem~\ref{thm:higher-central-extension}, applied with $s=1$, gives
\[
   \expn\widetilde G'
   \mid
     (\expn G)^{\mu_p(c)}
\]
with $c\geq 1$. 
By Proposition~\ref{prop:schur-cover-facts}(4),
we have $\widetilde G'\cong G\wed G$.
Since $\M(G)\le\widetilde G'$, the same bound
holds for $G\wed G$ and $\M(G)$. 
Since $p$ is odd, Blyth--Fumagalli--Morigi
\cite[Theorem~1.3(iii)]{BlythFumagalliMorigi2010}
give
\[
   G\otimes G\cong\nabla(G)\times(G\wed G),
   \qquad
   \nabla(G)\cong\nabla(G^{\mathrm{ab}}).
\]
Moreover, $
   \expn\nabla(G)\mid\expn G.
$ Hence the preceding bound for $G\wed G$ gives the same bound for $G\otimes G$.
\end{proof}

\begin{remark}\label{rem:quantitative-comparison}

Let $p$ be an odd prime and let $G$ be a finite $p$-group of
nilpotency class $c$. 
Bastos et al.~\cite[Theorem~1.4]{BastosMeloGoncalvesMonetta2022},
Sambonet~\cite[Theorem~1.1]{Sambonet2017},
and Antony--Komma--Thomas~\cite[Theorem~4.2]{AntonyKommaThomas2022}
obtained bounds for $G\otimes G$, $\M(G)$,
and $G\wed G$, respectively.
The last bound applies when $c\ge p$.
Theorem~\ref{thm:exterior-tensor-bounds} gives a common
bound for all three objects. The comparison is
summarized as follows.

\begin{table}[ht]
\centering
\caption{Exponent bounds for Schur multipliers,
nonabelian exterior squares, and tensor squares}
\label{tab:comparison-class-bounds}
\renewcommand{\arraystretch}{1.5}
\begin{tabularx}{\textwidth}{@{}l >{\centering\arraybackslash}X@{}}
\toprule
 Source & Bound \\
\midrule

BMGM
\cite[Theorem~1.4]{BastosMeloGoncalvesMonetta2022}
&
\(
   \expn(G\otimes G)
   \mid
   (\expn G)^{\lceil\log_p(c+1)\rceil}
\)
\\

\addlinespace

Sambonet
\cite[Theorem~1.1]{Sambonet2017}
&
\(
   \expn M(G)
   \mid
   (\expn G)^{\lfloor\log_{p-1}c\rfloor+1}
\)
\\

\addlinespace

AKT
\cite[Theorem~4.2]{AntonyKommaThomas2022}, $c\ge p$
&
\(
   \expn(G\wedge G)
   \mid
   (\expn G)^{
      1+\left\lceil
      \log_{p-1}\frac{c+1}{p+1}
      \right\rceil
   }
\)
\\

\addlinespace

Theorem~\ref{thm:exterior-tensor-bounds}
&
\(
   \expn X
   \mid
   (\expn G)^{
      1+\left\lceil
      \log_p\frac{c+1}{2p-1}
      \right\rceil
   }
   \text{ with } X \in\{G\wedge G,M(G),G\otimes G\}
\)
\\

\bottomrule
\end{tabularx}
\end{table}
\end{remark}

\subsection{Bogomolov multipliers}\label{sec:bogomolov}
Moravec's conjecture \cite[Conjecture~1.3]{Moravec2019Bogomolov} predicts $\expn\B(G)\mid\expn G$ for every finite group.   We prove that the Bogomolov multiplier exponent conjecture holds for finite \(p\)-groups of nilpotency class at most \(2p-1\). We also investigate the case where the nilpotency class exceeds \(2p-1\). Our results are as follows.
\begin{theorem}\label{thm:bog-quantitative}
Let $p$ be an odd prime and $G$ be a finite $p$-group. Write $q=\expn G$ and $c=\cl(G)$. Then, for $Y\in\{G\curlyw G,\B(G)\}$,
\[
   \expn Y\mid
   \begin{cases}
      q, & c\le 2p-1,\\[1mm]
      q^{\,1+\left\lceil\log_p\!\left(\dfrac{c}{2p-1}\right)\right\rceil},
      & c\ge 2p.
   \end{cases}
\]
\end{theorem}

\begin{proof}
The case $G=1$ is trivial, so assume $G\ne1$. Let $G^*$ be a CP cover of $G$. By Proposition~\ref{prop:cp-package}(3) and (4), we have
$(G^*)'\cong G\curlyw G$ and $G^*/Z(G^*)\cong G/Z(G)$. Thus $\expn(G^*/Z(G^*))\mid q$ and $\cl(G^*/Z(G^*))\le c-1$. Theorem~\ref{thm:higher-central-logarithmic}, applied to $G^*$ with $s=1$, therefore gives
\[
 \expn(G\curlyw G)=\expn(G^*)'\mid
   \begin{cases}
     q, & c\le 2p-1,\\[1mm]
      q^{\,1+\left\lceil\log_p\!\left(\dfrac{c}{2p-1}\right)\right\rceil},
      & c\ge 2p.
   \end{cases}
\]
Since $\B(G)\le G\curlyw G$, the same bound holds for $\B(G)$.
\end{proof}

\begin{remark}\label{rem:bogomolov-comparison}
Moravec \cite[Theorem~3.14]{Moravec2019Bogomolov} proved that
\(\expn\B(G)\mid\expn G\) for finite groups of nilpotency class at most
\(5\).  For finite \(p\)-groups with \(p\) odd, 
Theorem~\ref{thm:bog-quantitative} gives
an exponent bound for arbitrary nilpotency class.
The comparison is summarized as follows.

\begin{table}[ht]
\centering
\caption{Exponent bounds for the Bogomolov multiplier}
\label{tab:comparison-bogomolov-bounds}
\renewcommand{\arraystretch}{1.5}
\footnotesize
\begin{tabularx}{\textwidth}{@{}l c >{\centering\arraybackslash}X@{}}
\toprule
Source & Assumption & Bound \\
\midrule

Moravec
\cite[Theorem~3.14]{Moravec2019Bogomolov}
&
finite \(G\), \(c\le5\)
&
\(
   \expn\B(G)\mid\expn G
\)
\\

\addlinespace

Theorem~\ref{thm:bog-quantitative}
&
finite \(p\)-group, \(p\) odd, \(c\le2p-1\)
&
\(
   \expn\B(G)\mid\expn G
\)
\\

\addlinespace

Theorem~\ref{thm:bog-quantitative}
&
finite \(p\)-group, \(p\) odd, \(c\ge2p\)
&
\(
   \expn\B(G)
   \mid
   (\expn G)^{
      1+\left\lceil
      \log_p\dfrac{c}{2p-1}
      \right\rceil
   }
\)
\\

\bottomrule
\end{tabularx}
\end{table}
\end{remark}

\section{The critical class 2p-1}\label{sec:critical-class}

Throughout this section, $p$ is an odd prime and $G$ is a finite $p$-group with $\cl(G)\le2p-1$. We construct a homomorphism $\beta$ on a subgroup of $(G/G')\wed(G/G')$ such that $\beta=0$ if and only if $\expn\M(G)\mid\expn G$. When $d(G)\le3$, we show that $\beta=0$. When $d(G)=4$ and $\expn G=p$, we show that $|\operatorname{Im}\beta|\le p$.

\subsection{A degree-two power map}\label{subsec:critical-power-map}

Let $q=\expn G$ and choose a Schur cover
\begin{equation}\label{eq:critical-cover}
   1\longrightarrow \M(G)\longrightarrow \widetilde G\longrightarrow G\longrightarrow1,
   \qquad
   \M(G)\le Z(\widetilde G)\cap \widetilde G'.
\end{equation}
Since $\widetilde G/\M(G)\cong G$, we have $\g_{2p}(\widetilde G)\le \M(G)\le Z(\widetilde G)$ and $\widetilde G^q\le \M(G)$. Hence $\cl(\widetilde G)\le2p$. Let $B=G/G'$. Since $\M(G)\le \widetilde G'$, we have $G'\cong (\widetilde G/\M(G))'=\widetilde G'/\M(G)$, and therefore the projection $\widetilde G\to G$ induces an isomorphism $\widetilde G/\widetilde G'\cong G/G'=B$.

By \cite[Proposition~1]{Ellis1987Exterior},
we have
$
   B\wedge B\cong (G\wed G)/(K_1K_2),
$
where \(K_1\) and \(K_2\) are the images of \(G'\wedge G\) and \(G\wedge G'\), respectively. 
By Proposition~\ref{prop:schur-cover-facts}(4),
we have $G\wed G\cong\widetilde G'$.
Under this identification, the images of both \(G'\wedge G\) and \(G\wedge G'\) correspond to \(
   [\widetilde G',\widetilde G]
   =\g_3(\widetilde G).
\) Consequently, the natural maps induce an isomorphism
\begin{equation}\label{eq:critical-wedge}
   \theta:B\wedge B
   \xrightarrow{\ \cong\ }
   \widetilde G'/\g_3(\widetilde G).
\end{equation}
Via this isomorphism, we identify $B\wedge B$ with $\widetilde G'/\g_3(\widetilde G)$. Thus, if $a,b\in B$ and $x,y\in\widetilde G$ map to $a,b$, respectively, then $\theta(a\wedge b)=[x,y]\g_3(\widetilde G)$. 

Recall that an element of $B\wedge B$ is called \emph{decomposable} if it is of the form $a\wedge b$ for some $a,b\in B$. Thus a decomposable element $a\wedge b$ corresponds precisely to a degree-two class represented by a single commutator $[x,y]$ modulo $\g_3(\widetilde G)$. Define
\begin{equation}\label{eq:critical-W}
   W=\theta^{-1}(\frac{\M(G)\g_3(\widetilde G)}{\g_3(\widetilde G)})
   \le B\wedge B.
\end{equation}
Then $W\cong \M(G)/(\M(G)\cap\g_3(\widetilde G))$.
For $w\in W$, choose $m\in\M(G)$ such that
$\theta(w)=m\g_3(\widetilde G)$, and define
$\beta(w)=m^q$.

\begin{proposition}\label{prop:critical-obstruction}
The map $\beta:W\to\M(G)^q$ is a surjective homomorphism. Every decomposable
element of $W$ lies in $\ker\beta$. Moreover, $\expn\M(G)\mid q$ if and only if $\beta=0$.
\end{proposition}

\begin{proof}
We first note that
\begin{equation}\label{eq:critical-gamma3}
   \expn\g_3(\widetilde G)\mid q.
\end{equation}
Indeed, Theorem~\ref{thm:higher-central}, applied to $\widetilde G$ with $s=2$, gives $\expn\g_3(\widetilde G)\mid\expn(\widetilde G/Z_2(\widetilde G))\mid q$, since $\widetilde G^q\le \M(G)\le Z(\widetilde G)\le Z_2(\widetilde G)$.

For $w\in W$, the definition of $W$ ensures that $w$ is represented by $m\g_3(\widetilde G)$ for some $m\in \M(G)$.  If it is also represented by $m'\g_3(\widetilde G)$ with $m'\in \M(G)$, then $mm'^{-1}\in \M(G)\cap\g_3(\widetilde G)$, and hence $(mm'^{-1})^q=1$ by \eqref{eq:critical-gamma3}.  Thus $\beta$ is well defined.  Since $\M(G)$ is central, it is abelian, so $\beta$ is a homomorphism, and its image is exactly $\M(G)^q$.  Hence $\beta=0$ if and only if $\M(G)^q=1$.

It remains to consider decomposable elements.  Let $w=a\wedge b\in W$, with $a,b\in B$, and choose $x,y\in \widetilde G$ whose images in $B$ are $a,b$, respectively. Under the identification \eqref{eq:critical-wedge},  $\theta(w)=m\g_3(\widetilde G)$ for some $m\in \M(G)$. Then $\theta(w)=[x,y]\g_3(\widetilde G)=m\g_3(\widetilde G).$ Hence
$
   [x,y]=mh
$
for some $h\in\g_3(\widetilde G)$.  As $\M(G)$ is central and $h^q=1$ by \eqref{eq:critical-gamma3},
\begin{equation}\label{eq:critical-beta-commutator}
   \beta(w)=m^q=m^qh^q=(mh)^q=[x,y]^q.
\end{equation}

Let $K=\langle x,y\rangle$.
Since $[x,y]=mh\in\M(G)\g_3(\widetilde G)$
and $\M(G)\g_3(\widetilde G)
\trianglelefteq\widetilde G$, we have
$K'\le\M(G)\g_3(\widetilde G)$.
As $\M(G)\le Z(\widetilde G)$,
it follows inductively that
$
   \g_i(K)\le\g_{i+1}(\widetilde G), (i\ge3).
$
Thus $\cl(K)\le2p-1$.  Moreover, $K^q\le \widetilde G^q\cap K\leq\M(G)\cap K\le Z(K)$.  Theorem~\ref{thm:higher-central-extension} therefore gives $\expn K'\mid\expn(K/Z(K))\mid q$. In particular, $[x,y]^q=1$.  Equation \eqref{eq:critical-beta-commutator} gives $\beta(w)=1$, as required.
\end{proof}

\subsection{At most three generators}\label{subsec:critical-three}

\begin{proof}[Proof of Theorem~\ref{thm:critical-three-generators}(1)]

Let $q=\expn G$ and $B=G/G'$.  If $q=1$, then $G=1$ and there is nothing to prove.  Hence assume $q>1$.  Since $G$ is a finite $p$-group, $d(B)=d(G)\le3$ and $\expn B\mid q$. 
Choose generators $b_1,b_2,b_3$ of $B$ and a
surjective homomorphism
$\varphi:E=(\mathbb Z/q\mathbb Z)^3\to B$
sending $e_i$ to $b_i$, where $e_1,e_2,e_3$
are the generators of $E$. 
By functoriality of the exterior square,
$\varphi$ induces a surjective homomorphism
$E\wedge E\to B\wedge B$. Moreover $B\wedge B$ is generated by $b_i\wedge b_j$ $(1\le i<j\le3)$, and each such element is the image of $e_i\wedge e_j$. 

We show that every element of $E\wedge E$ is decomposable.  Write
\[
   \alpha
   =r\,e_1\wedge e_2
    +s\,e_1\wedge e_3
    +t\,e_2\wedge e_3,
   \qquad
   r,s,t\in\mathbb Z/q\mathbb Z.
\]
There is nothing to prove if $\alpha=0$.  Otherwise, let $u$ be maximal such that $r,s,t\in p^u(\mathbb Z/q\mathbb Z)$.  After permuting the basis, we may assume that $r=p^ur_0$, $s=p^us_0$, and $t=p^ut_0$, where $r_0$ is a unit modulo $q$. Then
\[
\begin{aligned}
   \alpha
   &=p^ur_0\,e_1\wedge e_2+p^us_0\,e_1\wedge e_3-p^ut_0\,e_3\wedge e_2\\
   &=e_1\wedge\bigl(p^ur_0e_2+p^us_0e_3\bigr)
     -t_0r_0^{-1}e_3\wedge\bigl(p^ur_0e_2+p^us_0e_3\bigr)\\
   &=\bigl(e_1-t_0r_0^{-1}e_3\bigr)
     \wedge\bigl(p^ur_0e_2+p^us_0e_3\bigr),
\end{aligned}
\]
where in the second equality we used \(e_3\wedge e_3=0\). Thus every element of $E\wedge E$, and hence every element of $B\wedge B$, is decomposable.  In particular, every element of the subgroup $W$ in \eqref{eq:critical-W} is decomposable.  Proposition~\ref{prop:critical-obstruction} gives $\expn\M(G)\mid q$.
\end{proof}

\subsection{Four generators of exponent p}\label{subsec:critical-four}

In this subsection, we investigate the four-generator case for  finite $p$-group $G$ of exponent $p$. We show that  $\M(G)$ have at most one cyclic factor of order $p^2$. We first prove some useful lemmas.

\begin{lemma}\label{lem:critical-derived-power}
Assume that $\expn G=p$ and $\cl(G)\le2p-1$, and retain the notation of  Schur cover $\widetilde G\to G$ from \eqref{eq:critical-cover}. Then
\begin{equation*}
   (\widetilde G')^p\le\g_{2p}(\widetilde G).
\end{equation*}
\end{lemma}

\begin{proof}
Set \(
   \overline{\widetilde G}
   =\widetilde G/\g_{2p}(\widetilde G).
\) Then $\cl(\overline{\widetilde G})\le2p-1$. Since $\widetilde G/\M(G)\cong G$ has exponent $p$, we have $\widetilde G^p\le\M(G)\le Z(\widetilde G)$. Hence \(
   \expn\bigl(
      \overline{\widetilde G}/Z(\overline{\widetilde G})
   \bigr)\mid p.
\) By Theorem~\ref{thm:higher-central-extension}, $
   \expn\overline{\widetilde G}'\mid \expn\bigl(
      \overline{\widetilde G}/Z(\overline{\widetilde G})
   \bigr)\mid p.
$ Since $\g_{2p}(\widetilde G)\le\widetilde G'$, we have \(
   \overline{\widetilde G}'
   =\widetilde G'/\g_{2p}(\widetilde G),
\) and therefore $
   (\widetilde G')^p\le\g_{2p}(\widetilde G),
$ as required.
\end{proof}

For the rest of this subsection, we assume that
\begin{equation*}
   \expn G=p,
   \qquad
   \cl(G)\le2p-1,
   \qquad
   d(G)=4.
\end{equation*}

\begin{corollary}\label{cor:critical-rank}
Let $p$ be an odd prime and $G$ a finite $p$-group
with $\expn G=p$, $\cl(G)\le2p-1$, and $d(G)=4$.
Then
\[
   \M(G)\cong C_{p^2}^{\,b}\times C_p^{\,a},
\]
for some $a,b\ge0$, where $b=\operatorname{rank}\beta$.
\end{corollary}

\begin{proof}
By Theorem~\ref{thm:exterior-tensor-bounds},
we have $\expn\M(G)\mid p^2$. Hence
$\M(G)\cong C_{p^2}^{\,b}\times C_p^{\,a}$
for some $a,b\ge0$.
By Proposition~\ref{prop:critical-obstruction},
$\operatorname{Im}\beta=\M(G)^p$.
Since $\beta$ is $\F_p$-linear, we have
$
   b=\dim_{\F_p}\M(G)^p
    =\operatorname{rank}\beta.
$
\end{proof}
Thus, to prove Theorem~\ref{thm:critical-three-generators}(2),
it remains to show that $\operatorname{rank}\beta\le1$.
Retain the notation of Subsection~\ref{subsec:critical-power-map}.  Since $G^p=1$, we have $G^pG'=G'$. By Burnside's basis theorem,  \(
   \dim_{\F_p}B=d(G)=4.
\) We therefore regard $B$ as a four-dimensional $\F_p$-vector space and write \[
   \Lambda^2_{\F_p}B
   =
   \frac{B\otimes_{\F_p}B}
        {\langle\,v\otimes v : v\in B\,\rangle}
\] for its second exterior power. We write $u\wedge v$ for the image of $u\otimes v$ in $\Lambda^2_{\F_p}B$. Since $B$ is elementary abelian, its group exterior square is canonically isomorphic to the usual exterior square of the underlying $\F_p$-vector space \cite{Ellis1987Exterior}. Thus  \(
   B\wedge B\cong\Lambda^2_{\F_p}B.
\) In particular, if $e_1,e_2,e_3,e_4$ is a basis of $B$, then
\[
   \Lambda^2_{\F_p}B
   =
   \left\langle
      e_i\wedge e_j : 1\le i<j\le4
   \right\rangle_{\F_p},
\]
and hence \(
   \dim_{\F_p}\Lambda^2_{\F_p}B=6.
\)  Write $T=\g_{2p}(\widetilde G)$. By \eqref{eq:critical-gamma3}, the group $T$ is an elementary abelian $p$-group, and hence an $\F_p$-vector space. Let $K=\F_{p^2}$, and define
\[
   B_K:=B\otimes_{\F_p}K,
   \qquad
   W_K:=W\otimes_{\F_p}K,
   \quad \text{ and }\quad
   T_K:=T\otimes_{\F_p}K.
\]
Let $\theta_K$ denote the $K$-linear extension of $\theta$ under the natural identification \(
    (B\wedge B)\otimes_{\F_p}K
    \cong\Lambda_K^2B_K.
\) We regard $W_K$ as a subspace of $\Lambda_K^2B_K$ and identify $B$, $W$, and $T$ with their images in $B_K$, $W_K$ and $T_K$ under the natural embeddings $x\mapsto x\otimes1$.

\begin{lemma}\label{lem:critical-quadratic-direction}
If $\operatorname{rank}\beta\ge2$,
then there exist $w_1,w_2\in W$,
$\lambda\in K\setminus\F_p$, and $u,v\in B_K$
such that $\beta(w_1)$ and $\beta(w_2)$
are linearly independent and
\[
   \lambda w_1+w_2=u\wedge v.
\]
\end{lemma}

\begin{proof}
Fix a basis $e_1,e_2,e_3,e_4$ of $B$. For \(
   z=\sum_{i<j}a_{ij}e_i\wedge e_j\in\Lambda_K^2B_K,
\) write
\[
   Q(z)=a_{12}a_{34}-a_{13}a_{24}+a_{14}a_{23}.
\]
A direct computation gives \(
z\wedge z
=2Q(z)\,e_1\wedge e_2\wedge e_3\wedge e_4.
\) 
Since $p$ is odd,
\cite[Lecture~6, Exercise~6.2]{Harris1992}
gives
\[
Q(z)=0
\quad\Longleftrightarrow\quad
z\wedge z=0
\quad\Longleftrightarrow\quad
z\text{ is decomposable}.
\]

Since $\operatorname{rank}\beta\ge2,$ we can choose $w_1,w_2\in W$ such that $t_1:=\beta(w_1)$ and $t_2:=\beta(w_2)$ are linearly independent in $\M(G)^p$. If $xw_1+yw_2\in\ker\beta$ with $x,y\in\F_p$, then $xt_1+yt_2=0$, so $x=y=0$. Thus $\beta$ is injective on $\langle w_1,w_2\rangle_{\F_p}$, and in particular $w_1,w_2$ are linearly independent. 
Proposition~\ref{prop:critical-obstruction}
implies that
\[
   \langle w_1,w_2\rangle_{\F_p}
   \cap\{u\wedge v:u,v\in B\}=\{0\}.
\]
Hence
$
   Q(xw_1+yw_2)\ne0
$
for all $(x,y)\in\F_p^2\setminus\{(0,0)\}$. In particular, $Q(w_1)\ne0$ and $Q(w_2)\ne0$. Since $Q$ is quadratic,
\[
   Q(xw_1+yw_2)
   =
   x^2Q(w_1)
   +xy\bigl(Q(w_1+w_2)-Q(w_1)-Q(w_2)\bigr)
   +y^2Q(w_2).
\]
Define \(
   f(Z):=Q(Zw_1+w_2).
\) 
The polynomial $f$ has degree two, and
$
   f(z)=Q(zw_1+w_2)\ne0
$
for all $z\in\F_p$. Thus $f$ is irreducible over $\F_p$.
Thus $f$ has no root in $\F_p$ and is
therefore irreducible over $\F_p$.

Let $\lambda\in K\setminus\F_p$ be a root of $f$.
The vectors $w_1,w_2$ remain linearly independent
over $K$, so $\lambda w_1+w_2\ne0$, while
$
   Q(\lambda w_1+w_2)=f(\lambda)=0.
$
Thus $\lambda w_1+w_2$ is decomposable, and hence there exist $u,v\in B_K$ such that
\begin{equation}\label{eq:critical-alpha}
   \lambda w_1+w_2=u\wedge v.
\end{equation}
\end{proof}

\begin{remark}\label{remark:motivation1}
By \cite[Chapter I, \S 2.2, Corollary 2]{Serre1973},
every three-dimensional subspace of $W$
contains a nonzero decomposable element.
Since every decomposable element of $W$
lies in $\ker\beta$, it follows that
$\operatorname{rank}\beta\le2$.

To obtain the sharper bound, suppose that
$\operatorname{rank}\beta=2$.
In the proof of Lemma~\ref{lem:critical-quadratic-direction},
we constructed a two-dimensional subspace
$\langle w_1,w_2\rangle_{\F_p}$
on which $\beta$ is injective.
Thus this subspace contains no nonzero decomposable
elements, whereas $\langle w_1,w_2\rangle_K$
contains the nonzero decomposable element
$\alpha=\lambda w_1+w_2$.

We will construct a map $F_K$ whose restriction
to $\langle w_1,w_2\rangle_K$ has trivial kernel,
and show that $F_K(\alpha)=0$.
Since $\alpha\ne0$, this gives a contradiction
and proves $\operatorname{rank}\beta\le1$.
\end{remark}

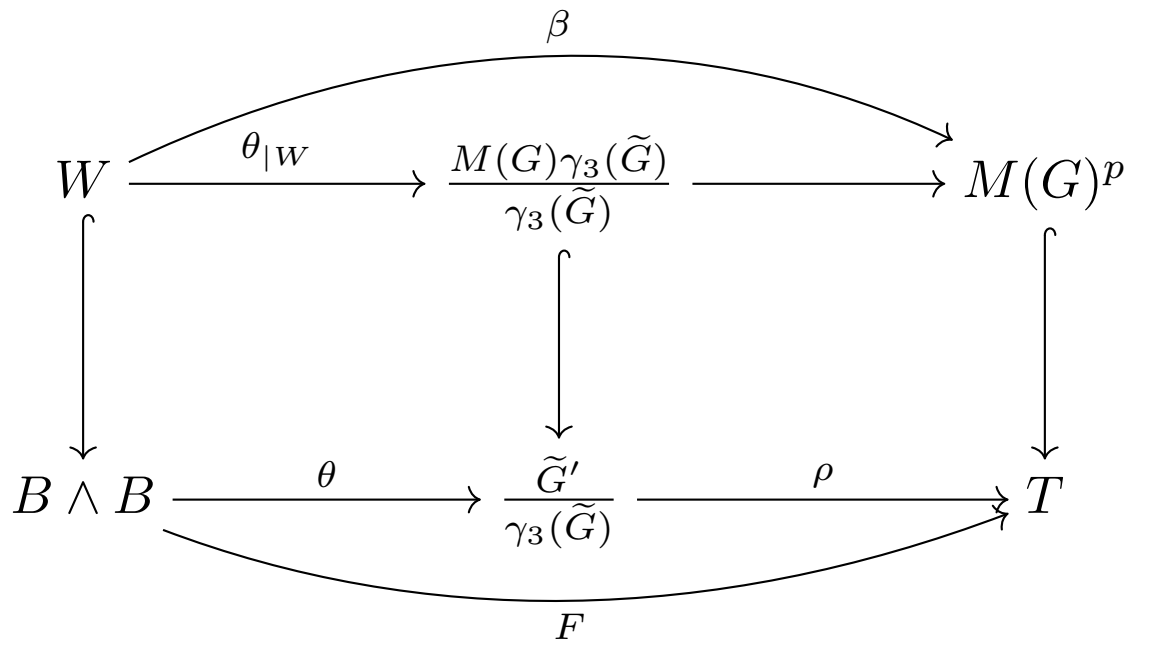
\begin{figure}[ht]
\centering
\begin{tikzcd}[ampersand replacement=\&, cramped]
    W \&\& {\frac{M(G)\gamma_3(\widetilde{G})}{\gamma_3(\widetilde{G})}} \&\& {M(G)^p} \\
    \\
    {B\wedge B } \&\& { \frac{\widetilde G'}{\gamma_3(\widetilde{G})}} \&\& T
    \arrow["{\theta_{\mid W}}", from=1-1, to=1-3]
    \arrow["\beta", curve={height=-30pt}, from=1-1, to=1-5]
    \arrow[hook, from=1-1, to=3-1]
    \arrow[from=1-3, to=1-5]
    \arrow[hook, from=1-3, to=3-3]
    \arrow[hook, from=1-5, to=3-5]
    \arrow["\theta", from=3-1, to=3-3]
    \arrow["F"', curve={height=24pt}, from=3-1, to=3-5]
    \arrow["\rho", from=3-3, to=3-5]
\end{tikzcd}
\caption{The critical obstruction diagram}
\label{fig:critical-obstruction-diagram}
\end{figure}

We now introduce the $p$-th-power map that will be used in the scalar-extension argument. By Lemma~\ref{lem:critical-derived-power} and the observation after  \eqref{eq:critical-cover},
\[
     (\widetilde G')^p\le T\le\M(G)\le Z(\widetilde G).
\]
Moreover $T\le\g_3(\widetilde G)$. By \eqref{eq:critical-gamma3}, we have \(
   \expn T\mid p.
\) Thus $T$ is naturally an $\F_p$-vector space. Define the $p$-th-power map
\[
   \rho:
   \widetilde G'/\g_3(\widetilde G)
   \longrightarrow T,
   \qquad
   \rho\bigl(c\g_3(\widetilde G)\bigr)=c^p.
\]
By Lemmas~\ref{lem:critical-derived-power}
and~\ref{lem:power-cosets}, the map $\rho$
is well defined. Composing $\rho$ with the
isomorphism $\theta$ in \eqref{eq:critical-wedge},
we obtain
\[
   F=\rho\circ\theta:
   B\wedge B\longrightarrow T.
\]
These maps are illustrated in
Figure~\ref{fig:critical-obstruction-diagram}.

\begin{lemma}\label{lem:critical-power-polynomial}

Let $w_1,w_2$ be as in
Lemma~\ref{lem:critical-quadratic-direction},
and write $t_i=\beta(w_i)$ for $i=1,2$, after choosing $\F_p$-bases of $B\wedge B$ and $T$, the coordinate functions of \(
   F:B\wedge B\longrightarrow T
\) are homogeneous polynomials of degree $p$. Extending these polynomials from $\F_p$ to $K$ defines a map \(
   F_K:\Lambda_K^2B_K\longrightarrow T_K
\) such that \(
   F_K(Xw_1+Yw_2)
   =
   X^p t_1+Y^p t_2
\text{ with }  X,Y \in K.
\)
\end{lemma}

\begin{proof}
Since $\dim_{\F_p}B=4$, we have $\dim_{\F_p}(B\wedge B)=6$. Extend $w_1,w_2$ to an $\F_p$-basis \(
   a_1=w_1,
   a_2=w_2,
   a_3,\ldots,a_6
\) of $B\wedge B$. Choose $c_i\in\widetilde G'$ such that \(
   \theta(a_i)=c_i\g_3(\widetilde G)
   \text{ with } 1\le i\le6.
\) Since $a_1=w_1$ and $a_2=w_2$ lie in $W$, we may choose $c_1,c_2\in\M(G)$. Hence, by the definition of $\beta$, \(
   c_i^p=t_i \text{ with }
i=1,2.
\) For $x=\sum_{i=1}^6x_i a_i\in B\wedge B$, choose integer representatives $0\le x_i<p$ and define \(
   c(x):=c_1^{x_1}\cdots c_6^{x_6}\in\widetilde G'.
\)  By the definitions of $\theta$ and $\rho$, we have
\[
F\left(\sum_{i=1}^6x_i a_i\right)=c(x)^p.
\]
We first work on the group side. By Proposition~\ref{prop:hall-collection}(1),
\[
   c(x)^p
   =
   \prod_{i=1}^6(c_i^{x_i})^p
   \prod_{r\ge2}\prod_{h\in\mathcal H_r}
   h(c_1^{x_1},\ldots,c_6^{x_6})^{e_h(p)},
\]
where $\mathcal H_r$ denotes the set of basic Hall commutators of weight $r$, and $e_h(p)$ is the corresponding exponent in the Hall collection formula. Since $c_i\in\widetilde G'=\g_2(\widetilde G)$, every weight-$r$ commutator in the $c_i^{x_i}$ lies in $\g_{2r}(\widetilde G)$. If $2\le r<p$, then Proposition~\ref{prop:hall-collection}(2) gives $p\mid e_h(p)$. Moreover, \(
   h(c_1^{x_1},\ldots,c_6^{x_6})
   \in\g_{2r}(\widetilde G)
   \le\g_3(\widetilde G)
\) and $\expn\g_3(\widetilde G)\mid p$. Hence the corresponding term is trivial.  If $r>p$, then \(
   h(c_1^{x_1},\ldots,c_6^{x_6})
   \in\g_{2p+2}(\widetilde G)=1.
\)  Thus 
\[
   c(x)^p
   =
   \prod_{i=1}^6(c_i^p)^{x_i}
   \prod_{h\in\mathcal H_p}
   h(c_1^{x_1},\ldots,c_6^{x_6})^{e_h(p)}.
\]
For $h\in\mathcal H_p$, let $n_i(h)$ denote the number of occurrences of the $i$th variable in $h$. Then
\(
   n_1(h)+\cdots+n_6(h)=p.
\) Since $h(c_1,\ldots,c_6)\in\g_{2p}(\widetilde G)$ and $\g_{2p+1}(\widetilde G)=1$, Lemma~\ref{lem:multilinearity} gives
\[
   h(c_1^{x_1},\ldots,c_6^{x_6})
   =
   h(c_1,\ldots,c_6)^{
      x_1^{n_1(h)}\cdots x_6^{n_6(h)}
   }.
\]
Therefore
\begin{equation}\label{eq:critical-hall-multiplicative}
\begin{aligned}
   c(x)^p
   ={}&
   \prod_{i=1}^6(c_i^p)^{x_i}\cdot \prod_{h\in\mathcal H_p}
   h(c_1,\ldots,c_6)^{
      e_h(p)
      x_1^{n_1(h)}\cdots x_6^{n_6(h)}
   }.
\end{aligned}
\end{equation}
We now write $T=\g_{2p}(\widetilde G)$ additively. 
Since $x_i^p=x_i$ for $x_i\in\F_p$,  \eqref{eq:critical-hall-multiplicative} becomes
\begin{equation}\label{eq:critical-hall-polynomial}
   c(x)^p
   =
   \sum_{i=1}^6x_i^p c_i^p
   +
   \sum_{h\in\mathcal H_p}
   e_h(p)
   x_1^{n_1(h)}\cdots x_6^{n_6(h)}
   h(c_1,\ldots,c_6).
\end{equation}

Choose an $\F_p$-basis of $T$. Since $n_1(h)+\cdots+n_6(h)=p$ for every $h\in\mathcal H_p$,
replacing $x_1,\ldots,x_6$ in \eqref{eq:critical-hall-polynomial} by indeterminates
$X_1,\ldots,X_6$ gives homogeneous coordinate polynomials of degree $p$ representing $F$.
These polynomial representatives are unique. Indeed, suppose that a homogeneous polynomial $P\in\F_p[X_1,\ldots,X_n]$ of degree $p$ vanishes on $\F_p^n$. Evaluating at the  basis vectors shows that the coefficients of all $X_i^p$ are zero. Thus $\deg_{X_i}P<p$ for every $i$, and hence $P=0$ by \cite[Lemma~2.1]{Alon1999}. Thus the homogeneous polynomial representation of $F$ is uniquely determined and is compatible with changes of bases. Extending its coefficients from $\F_p$ to $K$ therefore defines a polynomial map
\[
   F_K:\Lambda_K^2B_K\longrightarrow T_K,
\]
where $a_1,\ldots,a_6$ are regarded as a $K$-basis of $\Lambda_K^2B_K$.
Since $a_1=w_1$ and $a_2=w_2$, substituting $(X_1,\ldots,X_6)=(X,Y,0,\ldots,0)$ into
\eqref{eq:critical-hall-polynomial} gives
\[
\begin{aligned}
F_K(Xw_1+Yw_2)
={}&X^pc_1^p+Y^pc_2^p
+\sum_{h\in\mathcal H_p}
e_h(p)\left(\prod_{i=1}^6X_i^{n_i(h)}\right)
h(c_1,\ldots,c_6),
\end{aligned}
\]
where $(X_1,\ldots,X_6)=(X,Y,0,\ldots,0)$. 
Since $X_i=0$ for $3\le i\le6$, the product
$
   \prod_{i=1}^6 X_i^{n_i(h)}
$
vanishes unless $h$ involves only $c_1,c_2$.
In the latter case, $h(c_1,\ldots,c_6)=0$
since $c_1,c_2\in\M(G)\le Z(\widetilde G)$. Thus all Hall commutator terms vanish, and
\[
F_K(Xw_1+Yw_2)
=X^pt_1+Y^pt_2.
\]
\end{proof}

\begin{remark}\label{remark:motivation2}
Since $t_1,t_2$ are linearly independent over $K$,
\[
F_K(Xw_1+Yw_2)=0
\iff X^pt_1+Y^pt_2=0
\iff X^p=Y^p=0
\iff X=Y=0.
\]
In particular, for the nonzero decomposable element
$
\alpha=\lambda w_1+w_2=u\wedge v
$
given by Lemma~\ref{lem:critical-quadratic-direction},
we have $F_K(\alpha)\ne0$.
To obtain a contradiction, it remains to show that
$F_K(u\wedge v)=0$.
\end{remark}

We now turn to the remaining case. We first establish the following lemma using the Hall--Petrescu collection formula.
\begin{lemma}\label{lem:hall-p-power}
Let $H$ be a finite $p$-group with $\expn\g_3(H)\mid p$. If $x\in H$ and $d\in\g_2(H)$, then
\[
   (xd)^p
   \equiv
   x^p d^p [d,{}_{p-1}x]
   \pmod{\g_{p+2}(H)}.
\]
\end{lemma}

\begin{proof}
By the Hall--Petrescu collection formula \cite[Lemma 2.7(a)]{Moravec2007Power},
\[
   (xd)^p
   \equiv
   x^p d^p
   \prod_{i=1}^{p-1}
   [d,{}_i x]^{\binom{p}{i+1}}
   \pmod{C_{p,1}(x,d)},
\]
Since $d\in\g_2(H)$, the definition of $C_{p,1}(x,d)$ gives
\[
C_{p,1}(x,d)
\le \g_{p+2}(H)
\langle z^p : z\in\g_3(H)\rangle
=\g_{p+2}(H),
\]
where the equality follows from $\expn\g_3(H)\mid p$. Moreover, for $1\le i<p-1$, we have $[d,{}_i x]\in\g_3(H)$ and $p\mid\binom{p}{i+1}$, so these factors are trivial. Since $\binom pp=1$, only $[d,{}_{p-1}x]$ remains. Thus
\[
   (xd)^p
   \equiv
   x^p d^p [d,{}_{p-1}x]
   \pmod{\g_{p+2}(H)}.
\]
\end{proof}

\begin{lemma}\label{lem:critical-balanced}
Let $H$ be a finite $p$-group with $H^p\le Z(H)$ and $\cl(H)\le2p$. Then there exist commutator words $c_1(X,Y),\ldots,c_N(X,Y)$ of multiweight $(p,p)$ and integers $\varepsilon_1,\ldots,\varepsilon_N$ such that
\[
   [x,y]^p
   =
   \prod_{\nu=1}^N c_\nu(x,y)^{\varepsilon_\nu}
   \qquad(x,y\in H).
\]
\end{lemma}

\begin{proof}
Apply the construction of Section~\ref{sec:induction} with $q=p$ and
$k=2p$. Denote the group in \eqref{eq:two-generator-V} by $V$, with
distinguished generators $\bar a,\bar b$. By Proposition~\ref{prop:universal-reduction}, we have $\bar a^p=\bar b^p=1$, $V^p\le Z(V)$, and $\cl(V)\le2p$. 
Since $V^p\le Z(V)\le Z_2(V)$, Theorem~\ref{thm:higher-central}, applied with $s=2$, gives $\expn\g_3(V)\mid p$.
By Lemma~\ref{lem:power-depth}, \(
z:=[\bar b,\bar a]^p\in\g_{2p}(V).
\) Since $\g_{2p}(V)\le\g_3(V)$ and $\cl(V)\le2p$, the group $\g_{2p}(V)$ is elementary abelian, and we write it additively.

For $r,s\in\F_p^\times$, let $\sigma_{r,s}$ be the automorphism of
$V$ induced by
\[
   \bar a\longmapsto\bar a^r
   \quad \text{and } \quad
   \bar b\longmapsto\bar b^s.
\]
Since $\g_{2p}(V)$ is characteristic in $V$, each $\sigma_{r,s}$
restricts to an $\F_p$-linear automorphism of $\g_{2p}(V)$. By
Proposition~\ref{prop:power-scaling},
$\sigma_{r,s}(z)=z^{rs}$. In additive notation, this becomes
$\sigma_{r,s}(z)=rs\,z$.
Since $V=\langle\bar a,\bar b\rangle$ and
$\g_{2p+1}(V)=1$, the group $\g_{2p}(V)$ is
generated by the weight-$2p$ commutator words
in $\bar a,\bar b$. We may therefore choose
an $\F_p$-basis $c_1,\ldots,c_N$ consisting
of such words and write \(
   z=\sum_{t=1}^N\lambda_t c_t.
\) Each basis $c_t$ of multiweight $(i_t,j_t)$ satisfies
\[
   \sigma_{r,s}(c_t)=r^{i_t}s^{j_t}c_t,
\]
whereas $\sigma_{r,s}(z)=rs\cdot z$. Comparing coefficients, every basis occurring with a nonzero coefficient must satisfy \(
   r^{i_t}s^{j_t}=rs.
\) Hence
\[
   i_t\equiv j_t\equiv1\pmod{p-1}.
\]
Since $i_t+j_t=2p$, the only possible multiweights are
\[
   (1,2p-1),\qquad(p,p) \quad \text{and}\quad(2p-1,1).
\]
We next consider the multiweight $(1,2p-1)$. Set $d=[\bar b,\bar a]$. Since
$\bar b^{\bar a}=\bar b d$ and $d^p=z\in\g_{2p}(V)\le\g_{p+2}(V)$, Lemma~\ref{lem:hall-p-power} gives
\[
\bar b^p
=(\bar b d)^p
\equiv
\bar b^p[d,{}_{p-1}\bar b]
\pmod{\g_{p+2}(V)}.
\]
Hence \(
[d,{}_{p-1}\bar b]\in\g_{p+2}(V),
\) and therefore \(
[d,{}_{2p-2}\bar b]
\in\g_{2p+1}(V)=1.
\)

Consider a commutator of multiweight $(1,2p-1)$ in $\bar a,\bar b$. If the unique occurrence of $\bar a$ is not among the first two entries, then the first two entries are both $\bar b$, so the commutator is trivial. Thus the only possibly nontrivial such commutators are
\[
[\bar b,\bar a,{}_{2p-2}\bar b]
\quad\text{or}\quad
[\bar a,\bar b,{}_{2p-2}\bar b].
\]
The first is trivial by the preceding identity. The second is its inverse modulo $\g_{2p+1}(V)=1$, and hence is also trivial. Also, all commutators of multiweight $(2p-1,1)$ are trivial. It follows that there exist commutator words $c_1(X,Y),\ldots,c_N(X,Y)$ of multiweight $(p,p)$ and integers $\varepsilon_1,\ldots,\varepsilon_N$ such that \(
z=[\bar b,\bar a]^p
=\prod_{\nu=1}^N
c_\nu(\bar a,\bar b)^{\varepsilon_\nu}.
\)
By Proposition~\ref{prop:universal-reduction}(2), the natural map $U\to V$ induces an isomorphism $U'\cong V'$. Applying Proposition~\ref{prop:universal-reduction}(1) to  $x,y\in H$ yields \(
[y,x]^p
=\prod_{\nu=1}^N
c_\nu(x,y)^{\varepsilon_\nu}.
\) Interchanging $x$ and $y$ gives the desired result.

\end{proof}

\begin{remark}
\label{remark mu}
 By the observation after  \eqref{eq:critical-cover}, we applies Lemma~\ref{lem:critical-balanced}   to $\widetilde G$. Thus there exist commutator words $c_1(X,Y),\ldots,c_N(X,Y)$ of multiweight $(p,p)$ and integers $\varepsilon_1,\ldots,\varepsilon_N$ such that
\begin{equation}\label{eq:critical-balanced-fixed}
   [x,y]^p
   =
   \prod_{\nu=1}^N c_\nu(x,y)^{\varepsilon_\nu}
   \qquad(x,y\in\widetilde G).
\end{equation}
Moreover, let $c(X,Y)$ be a commutator word of multiweight $(r,s)$. For $a_1,\ldots,a_r,b_1,\ldots,b_s$, we write \(
c(a_1,\ldots,a_r;b_1,\ldots,b_s)
\) for the word obtained from $c(X,Y)$ by replacing the successive occurrences of $X$ and $Y$ with $a_1,\ldots,a_r$ and $b_1,\ldots,b_s$, respectively, while preserving the original commutator structure.
\end{remark}

\begin{lemma}\label{lem:critical-balanced-extension}
With the notation of Lemma~\ref{lem:critical-power-polynomial}, there exists a map
\[
   \Theta_K:B_K\times B_K\longrightarrow T_K
\]
such that
\[
   F_K(r\wedge s)=\Theta_K(r,s)
 \text{ and }r,s\in B_K,
\] and the coordinate functions of $\Theta_K$ are bihomogeneous polynomials of bidegree $(p,p)$.
\end{lemma}

\begin{proof}
Let $   \pi:\widetilde G\longrightarrow G/G'=B$ be the natural map. For $r,s\in B$, choose $x,y\in\widetilde G$ such that \(
   \pi(x)=r \text{ and }\pi(y)=s.
\) By \eqref{eq:critical-wedge} and the definition of $F=\rho\circ\theta$,
\begin{equation}\label{eq:critical-decomposable-Fp}
   F(r\wedge s)=[x,y]^p\in T.
\end{equation}
Choose an $\F_p$-basis $e_1,\ldots,e_4$ of $B$, and choose $z_1,\ldots,z_4\in\widetilde G$ such that \(
   \pi(z_i)=e_i
   \text{ with } 1\le i\le4.
\) For \(
   R=(R_1,\ldots,R_4)\text{ and }
   S=(S_1,\ldots,S_4)\in\F_p^4,
\) choose integer representatives $0\le R_i,S_i<p$. 
Write
\[
\begin{alignedat}{2}
   r(R)&:=\sum_{i=1}^4 R_i e_i,
   \qquad &s(S)&:=\sum_{i=1}^4 S_i e_i,\\
   x(R)&:=z_1^{R_1}\cdots z_4^{R_4},
   &y(S)&:=z_1^{S_1}\cdots z_4^{S_4}.
\end{alignedat}
\]
Then $\pi(x(R))=r(R)$ and $\pi(y(S))=s(S)$.
Fix $1\le\nu\le N$. By the choice of $c_\nu$ in \eqref{eq:critical-balanced-fixed}, it has multiweight $(p,p)$. 
Since $\g_{2p+1}(\widetilde G)=1$,
Lemma~\ref{lem:multilinearity} gives
\[
\begin{aligned}
c_\nu(x(R),y(S))
={}&
\prod_{\substack{i_1,\ldots,i_p\in\{1,\ldots,4\}\\
                   j_1,\ldots,j_p\in\{1,\ldots,4\}}}
 c_\nu
 \bigl(
 z_{i_1},\ldots,z_{i_p};
 z_{j_1},\ldots,z_{j_p}
 \bigr)^{
 R_{i_1}\cdots R_{i_p}
 S_{j_1}\cdots S_{j_p}}.
\end{aligned}
\]
Since each commutator value lies in $T=\g_{2p}(\widetilde G)$, we write the product additively as
\[
c_\nu(x(R),y(S))
=
\sum_{\substack{i_1,\ldots,i_p\in\{1,\ldots,4\}\\
j_1,\ldots,j_p\in\{1,\ldots,4\}}}
R_{i_1}\cdots R_{i_p}S_{j_1}\cdots S_{j_p}\,
c_\nu(z_{i_1},\ldots,z_{i_p};z_{j_1},\ldots,z_{j_p}).
\]
 We extend the coefficients to $T_K=T\otimes_{\F_p}K.$  We define
\[
   C_{\nu,K}:B_K\times B_K\longrightarrow T_K
\]
by
\[
\begin{aligned}
C_{\nu,K}(r(R),s(S))
:=
\sum_{\substack{i_1,\ldots,i_p\in\{1,\ldots,4\}\\
j_1,\ldots,j_p\in\{1,\ldots,4\}}}
\quad R_{i_1}\cdots R_{i_p}S_{j_1}\cdots S_{j_p}\,
\Bigl(
c_\nu(z_{i_1},\ldots,z_{i_p};z_{j_1},\ldots,z_{j_p})
\otimes1
\Bigr),
\end{aligned}
\]
where $r(R),s(S)\in B_K$ are given by the same coordinate formulas as above. The preceding identity gives \(
   C_{\nu,K}(r(R),s(S))
   =c_\nu(x(R),y(S))\otimes1.
\) Writing $\overline{\varepsilon_\nu}\in\F_p$ for the image of $\varepsilon_\nu$. Define $\Theta_K:B_K\times B_K\to T_K$ by
\[
   \Theta_K(r,s)
   :=
   \sum_{\nu=1}^N
   \overline{\varepsilon_\nu}C_{\nu,K}(r,s).
\]
By construction, the coordinate functions of $\Theta_K$ are bihomogeneous polynomials of bidegree $(p,p)$. Let $\Theta$ denote the same polynomial expression over $\F_p$, it follows from \eqref{eq:critical-balanced-fixed} and \eqref{eq:critical-decomposable-Fp} that
\begin{equation}\label{eq:critical-balanced-Fp}
   F\bigl(r(R)\wedge s(S)\bigr)
   =
   \Theta(r(R),s(S)).
\end{equation}
We now show that \eqref{eq:critical-balanced-Fp} is an identity of polynomials. Write
\[
   r(R)\wedge s(S)
   =
   \sum_{j=1}^6q_j(R,S)a_j.
\]
Each $q_j(R,S)$ is bilinear in $R$ and $S$. By Lemma~\ref{lem:critical-power-polynomial}, the coordinate functions of $F$ are homogeneous polynomials of degree $p$. Hence the coordinate functions of \(
   F\bigl(r(R)\wedge s(S)\bigr)
\) are bihomogeneous polynomials of bidegree $(p,p)$. Since the coordinate functions of $C_{\nu,K}(r(R),s(S))$ are bihomogeneous polynomials of bidegree $(p,p)$, the same holds for $\Theta(r(R),s(S))$.
Consider the coordinate polynomials of
\[
H(R,S):=
F(r(R)\wedge s(S))-\Theta(r(R),s(S)).
\]
They are bihomogeneous of bidegree $(p,p)$ and vanish on $\F_p^4\times\F_p^4$. Fix $S_0\in\F_p^4$. Since $H(R,S_0)$ is homogeneous of degree $p$ in $R$ and vanishes for  $R\in\F_p^4$, the uniqueness established in the proof of Lemma~\ref{lem:critical-power-polynomial} implies that it is the zero polynomial in $R$. Thus every coefficient of $H$, regarded as a polynomial
in $R$, vanishes on $\F_p^4$. Each such coefficient is homogeneous of degree $p$ in $S$, so applying
the same uniqueness result again gives \(
   H=0.
\) Therefore \eqref{eq:critical-balanced-Fp} is an identity of polynomials over $\F_p$. Extending scalars to $K$ gives
\begin{equation}
\label{eq:critical-balanced-scalar}
F_K(r\wedge s)=\Theta_K(r,s).
\end{equation}
\end{proof}

\begin{lemma}\label{lem:critical-W-vanishing}
With the notation of Lemma~\ref{lem:critical-balanced-extension}, if  $u,v\in B_K$,
$\alpha=u\wedge v\in W_K$,  
then $
   F_K(\alpha)=0.
$
\end{lemma}

\begin{proof}
By \eqref{eq:critical-balanced-scalar}
and the definition of $\Theta_K$, we have
\[
   F_K(\alpha)
   =\Theta_K(r,s)
   =\sum_{\nu=1}^N
      \overline{\varepsilon_\nu}C_{\nu,K}(r,s).
\]
Thus it suffices to show that
$C_{\nu,K}(r,s)=0$ for all $1\le\nu\le N$.

Recall that each $C_{\nu,K}$ is obtained from
a commutator of weight $2p$
by extending the commutator maps on successive
lower-central quotients to $K$.
Since
\[
   [\g_4(\widetilde G),
   \underbrace{\widetilde G,\ldots,\widetilde G}_{2p-3}]
   \le \g_{2p+1}(\widetilde G)=1,
\]
it suffices to show that the corresponding
degree-three bracket vanishes in
\[
   \bigl(\g_3(\widetilde G)/\g_4(\widetilde G)\bigr)
   \otimes_{\F_p}K.
\]
Recall that $\theta_K$ is the $K$-linear extension
of the isomorphism $\theta$ in
\eqref{eq:critical-wedge}
(see Figure~\ref{fig:critical-obstruction-diagram}).
Fix $1\le\nu\le N$. If the first two entries
of $c_\nu(X,Y)$ are the same variable, then
$c_\nu$ is trivial. Otherwise, they are
$X,Y$ or $Y,X$, and the corresponding
degree-two bracket is
$\pm\theta_K(\alpha)$.
Thus it suffices to show that the brackets
of $\theta_K(\alpha)$ with $r$ and $s$
vanish in $
   \bigl(\g_3(\widetilde G)/\g_4(\widetilde G)\bigr)
   \otimes_{\F_p}K.
$

The group commutator induces an $\F_p$-bilinear map
\[
\begin{aligned}
   \Phi:
   \g_2(\widetilde G)/\g_3(\widetilde G)
   \times B
   &\longrightarrow
   \g_3(\widetilde G)/\g_4(\widetilde G),\\
   (c\g_3(\widetilde G),b)
   &\longmapsto
   [c,x]\g_4(\widetilde G),
\end{aligned}
\]
where $c\in\g_2(\widetilde G)$ and
$x\in\widetilde G$ represents $b\in B$.
Since $\M(G)\le Z(\widetilde G)$,
we have $\Phi(\theta(W),B)=0$.
By $\F_p$-bilinearity, extending scalars
to $K$ yields
\begin{equation}\label{eq:critical-Phi-vanishing}
   \Phi_K(\theta_K(w),b)=0
   \qquad(w\in W_K,\ b\in B_K).
\end{equation}
Since $\alpha\in W_K$, we obtain
\[
   \Phi_K(\theta_K(\alpha),r)
   =\Phi_K(\theta_K(\alpha),s)=0.
\]
Hence $C_{\nu,K}(r,s)=0$ for every $\nu$,
and therefore $F_K(\alpha)=0$.
\end{proof}

\begin{proof}[Proof of Theorem~\ref{thm:critical-three-generators}(2)]
Suppose that $\operatorname{rank}\beta\ge2$.
Choose $w_1,w_2,\lambda,u,v$ as in
Lemma~\ref{lem:critical-quadratic-direction}.
By Lemmas~\ref{lem:critical-W-vanishing}
and~\ref{lem:critical-power-polynomial}, we have
\begin{equation}\label{eq:critical-plane-contradiction}
   0=F_K(u\wedge v)
    =F_K(\lambda w_1+w_2)
    =\lambda^p\beta(w_1)+\beta(w_2),
\end{equation}
contradicting the $K$-linear independence of
$\beta(w_1)$ and $\beta(w_2)$.
Hence $\operatorname{rank}\beta\le1$,
and the conclusion follows from
Corollary~\ref{cor:critical-rank}.
\end{proof}

\section*{Declaration of AI use}

The authors used ChatGPT (GPT-6, OpenAI) to assist with
language editing, manuscript organization,
and exploratory discussions and preliminary
checks of mathematical arguments.
All mathematical results and proofs in the final
manuscript were independently verified by the
authors, who take full responsibility for the content.

\bibliographystyle{abbrv}
\bibliography{schur_exponent_references}

\end{document}